\documentclass{birkjour}
\usepackage{mathabx}
\usepackage{eufrak}
\usepackage{amsfonts}
\usepackage{amssymb}
\usepackage{amsmath}
\usepackage{amssymb}
\usepackage[dvipsnames]{xcolor}
 \usepackage{xcolor}
 \usepackage{empheq}
\usepackage[most]{tcolorbox}
 \newtheorem{thm}{Theorem}[section]
 \newtheorem{cor}[thm]{Corollary}
 \newtheorem{lem}[thm]{Lemma}
 \newtheorem{prop}[thm]{Proposition}
 \theoremstyle{definition}
 \newtheorem{defn}[thm]{Definition}
 \theoremstyle{remark}

 \numberwithin{equation}{section}
\newtcbox{\mymath}[1][]{%
    nobeforeafter, math upper, tcbox raise base,
    enhanced, colframe=blue!30!black,
    colback=blue!30, boxrule=1pt,
    #1}
\definecolor{fen}{rgb}{0.31, 0.47, 0.26}

\begin{document}

%
%
%
%
%
%
%
%
%

\title[\footnotesize{Topologizability and Power Boundedness of Convolution and Toeplitz Operators}]
 {{Topologizability and Power Boundedness of Convolutions and Toeplitz Operators on Power Series Spaces}}

\author{Nazl\i\;Do\u{g}an}
\address{Fatih Sultan Mehmet Vakıf University, 
    34445 Istanbul, Turkey}
\email{ndogan@fsm.edu.tr}



\dedicatory{Dedicated to the memory of Vyacheslav Pavlovich Zakharyuta}

\begin{abstract}
We characterize the topologizability and power boundedness of convolution and dual convolution operators on power series spaces. We determine necessary conditions for a Toeplitz operator to be m-topologizable, and power bounded on $\Lambda_{1}(n)$ and $\Lambda_{\infty}(n)$, and consequently on $H(\mathbb{C})$ and $H(\mathbb{D})$.
\end{abstract}

\maketitle

\section{Introduction}
Toeplitz operators have been thoroughly studied in the context of Hilbert spaces, such as the Hardy and Bergman spaces, primarily because they create a bridge between function theory and operator theory. In the Hardy space $H^2(\mathbb{D})$, a Toeplitz operator corresponds to a Toeplitz matrix. Furthermore, every bounded linear operator on $H^2(\mathbb{D})$ whose matrix representation is Toeplitz arises from a Toeplitz operator. In the last two decades, the study of Toeplitz operators — those for which the matrix representation is of Toeplitz — has been generalized to broader settings beyond Hilbert spaces. For instance, Domański and Jasiczak \cite{DJ} developed such a framework for $A(\mathbb{R})$, the space of real analytic functions on $\mathbb{R}$, which is a non-metrizable and non-normable topological vector space. Later, in \cite{J1}, \cite{J2}, Jasiczak extended this line of study to Fréchet spaces, specifically to the space of entire functions $H(\mathbb{C})$ and to the space of holomorphic functions on the unit disc $H(\mathbb{D})$. The space of entire functions $H(\mathbb{C})$ and the space of holomorphic functions on the unit disc $H(\mathbb{D})$ are isomorphic to power series space of infinite type $\Lambda_{\infty}(n)$ and of finite type $\Lambda_{1}(n)$. In \cite{N1}, we extended the definition of Toeplitz operators to a broader class of power series spaces, including both finite and infinite types. Given a Toeplitz operator $T_{\theta,\beta}$ on a power series space, it expressed as $T_{\theta,\beta} = \widehat{T}_{\theta} + \widecheck{T}_{\beta}$, where $\widehat{T}_{\theta}$ and $\widecheck{T}_{\beta}$ are convolution-type operators. As shown at the end of Section 3.1 and 3.2, based on the identities $\widehat{T}_{\theta}x = \theta * x$ and $\widecheck{T}_{\beta}x = \beta \star x$, we refer to $\widehat{T}_{\theta}$ as the convolution operator with symbol $\theta$, and $\widecheck{T}_{\beta}$ as the dual convolution operator with symbol $\beta$  (where the dual convolution $\star$ is defined at the end of Section 3.2).

The terminology of topologizable and m-topologizable operators was first formulated by \.{Z}elazko in [13], from the perspective of equipping the space of continuous linear operators $\mathcal{L}(E)$, where $E$ is a non-normable locally convex space, with a compatible locally convex topology. In recent years, topologizability, power boundedness, mean ergodicity, and related properties have been extensively explored for continuous linear operators on locally convex spaces, especially Fr\'echet spaces, see, for instance, \cite{ABR}, \cite{ABR2}, \cite{B}, \cite{BD}, \cite{N4} \cite{K1}, \cite{K2}.

This paper aims to study topologizability, power boundedness of convolution, dual convolution, and Toeplitz operators on power series spaces. The paper is organized as follows. In Section 2, we establish the notational framework and gather essential lemmas and propositions that support the main results. Section 3 is devoted to the definition of convolution, dual convolution, and Toeplitz operators on power series spaces. In Section 4, we show that convolution operators are always m-topologizable and characterize their power boundedness through necessary and sufficient conditions. In Section 5, the necessary and sufficient conditions are derived for the dual convolution operators to be topologizable, m-topologizable, and power bounded. In Section 6, we establish necessary conditions for the Toeplitz operator defined on $\Lambda_{1}(n)$ and $\Lambda_{\infty}(n)$ to be m-topologizable and power bounded. This, in turn, will allow us to identify necessary conditions for the operator to be topologizable,
m-topologizable, and power bounded on the spaces $H(\mathbb{C})$ and $H(\mathbb{\mathbb{D}})$.
\section{Preliminaries}
This section outlines the essential background required throughout the paper. We adhere to the terminology and notation of \cite{MV}.

Fréchet spaces are complete Hausdorff locally convex spaces whose topology can be described by a countable system of seminorms $(\|\cdot\|_{k})_{k\in \mathbb{N}}$. A grading on a Fr\'echet space $E$ is a sequence of seminorms $\{\|\cdot\|_{n}:n\in \mathbb{N}\}$ which are increasing, that is,
$$\|x\|_{1}\leq \|x\|_{2}\leq \|x\|_{3}\leq \dots$$
for each $x\in E$ and which defines the topology. Every Fr\'echet space admits a grading. A graded Fr\'echet space is a Fr\'echet space with a choice of grading. For more information, see \cite{H}. In this paper, unless stated otherwise, we will assume that all Fr\'echet spaces are graded Fr\'echet spaces. 

A matrix $(a_{n,k})_{k,n\in \mathbb{N}}$ of non-negative numbers is called a Köthe matrix if the following conditions hold:
\begin{itemize}
\item[1.] For each $n\in \mathbb{N}$ there exists a $k\in \mathbb{N}$ with $a_{n,k}>0$. 
\item[2.] $a_{n,k}\leq a_{n,k+1}$ for each $n,k\in \mathbb{N}$.
\end{itemize}
Given a Köthe matrix $(a_{n,k})_{n,k\in \mathbb{N}}$, the space
\[ K(a_{n,k})=\bigg\{ x=(x_{n})_{n\in \mathbb{N}}\; \bigg|\;\|x\|_{k}:=\sum^{\infty}_{n=1}|x_{n}|a_{n,k}< \infty \;\text{for all}\;k\in \mathbb{N} \bigg\}\]
is called a Köthe space. Every Köthe space is a Fr\'echet space given by the semi-norms
in its definition. According to Proposition 27.3 of \cite{MV}, the dual space of a Köthe space $K(a_{n,k})$ is
\[ (K(a_{n,k}))^{\prime}=\bigg\{ y=(y_{n})_{n\in \mathbb{N}}\; \bigg|\; \sup_{n\in \mathbb{N}}|y_{n}a_{n,k}^{-1}|<+\infty\; \text{for some}\; k\in \mathbb{N}\bigg\}.\]


Let $\alpha=\left(\alpha_{n}\right)_{n\in \mathbb{N}}$ be a non-negative increasing sequence with $\displaystyle \lim_{n\rightarrow \infty} \alpha_{n}=\infty$. The power series space of finite type associated with $\alpha$ is defined as
$$\Lambda_{1}\left(\alpha\right):=\left\{x=\left(x_{n}\right)_{n\in \mathbb{N}}: \;\left\|x\right\|_{k}:=\sum^{\infty}_{n=1}\left|x_{n}\right|e^{-{1\over k}\alpha_{n}}<\infty \textnormal{ for all } k\in \mathbb{N}\right\},$$
and the power series space of infinite type associated with $\alpha$ is defined as
$$\displaystyle \Lambda_{\infty}\left(\alpha\right):=\left\{x=\left(x_{n}\right)_{n\in \mathbb{N}}:\; \left\|x\right\|_{k}:=\sum^{\infty}_{n=1}\left|x_{n}\right|e^{k\alpha_{n}}<\infty \textnormal{ for all } k\in \mathbb{N}\right\}.$$
Power series spaces constitute an important class of Köthe spaces. They include the spaces of holomorphic functions on $\mathbb{C}^{d}$ and on the polydisc $\mathbb{D}^{d}$, where $\mathbb{D}$ denotes the unit disk in $\mathbb{C}$ and $d\in \mathbb{N}$:
\[\mathcal{H}(\mathbb{C}^{d})\cong \Lambda_{\infty}(n^{\frac{1}{d}}) \quad\quad\text{and} \quad\quad \mathcal{H}(\mathbb{D}^{d})\cong \Lambda_{1}(n^{\frac{1}{d}}).
\]
The nuclearity of a power series space of finite type $\Lambda_{1}(\alpha)$ and of infinite type $\Lambda_{\infty}(\alpha)$ is equivalent to conditions $\displaystyle\lim_{n\to \infty} \frac{\ln n}{\alpha_{n}}=0$ and $\displaystyle\sup_{n\in \mathbb{N}} \frac{\ln n}{\alpha_{n}}<+\infty$,
respectively.  A sequence $\alpha$ is called stable if
\begin{equation}\label{E3.5}
\sup_{n\in \mathbb{N}} {\alpha_{2n}\over \alpha_{n}}<\infty.
\end{equation}
We want to note that the sequence $(n^{\frac{1}{d}})_{n\in \mathbb{N}}$ is stable and $\lim_{n\to \infty}n^{-\frac{1}{d}} \ln n =0$ for all $d\in \mathbb{N}$. Hence, the associated exponent sequence of $\mathcal{H}(\mathbb{C}^{d})$ and $\mathcal{H}(\mathbb{D}^{d})$ is stable and these spaces are nuclear for every $d\in \mathbb{N}$. 
 
Let $E$ be a Fr\'echet space. A linear map $T:E\to E$ is called continuous 
if for every $k\in \mathbb{N}$ there exists $p\in \mathbb{N}$ 
and $C_k>0$ such that 
\[\|Tx\|_k\leq C_k\|x\|_p\]
for all $x\in E$.  $\mathcal{L}(E)$ denotes the vector space of all continuous linear maps from $E$ to $E$. The composition of an operator $T\in \mathcal{L}(E)$ with itself n-times is denoted by $T^{n}$. A linear map $T: E\to F$ is called compact
if $T(U)$ is precompact in $F$ where $U$ is a neighborhood of zero of E. 

Throughout this paper, $e_{n}$ stands for the sequence with 1 in the $n^{th}$ component and zeros elsewhere. In any Köthe space, the sequence $(e_{n})_{n\in \mathbb{N}}$ forms the canonical Schauder basis.

The following lemma, due to Crone and Robinson (see Lemma 1 in \cite{CR75}), is a useful tool for analyzing the continuity and compactness of operators between Köthe spaces.
\begin{lem}\label{Crone} 
Let $K(a_{n,k})$ and $K(b_{n,k})$ be Köthe spaces.
\begin{itemize}
\item[a.] $T: K(a_{n,k}) \to K(b_{n,k})$ is a linear continuous operator if and only if for each $k\in \mathbb{N}$ there exists $m\in \mathbb{N}$ such that
\[ \sup_{n\in \mathbb{N}} {\frac{\|Te_{n}\|_{k}}{\|e_{n}\|_{m}}<\infty }.\]
\item[b.] If $K(b_{n,k})$ is Montel, then $T: K(a_{n,k}) \to K(b_{n,k})$ is a compact operator if and only if there exists $m\in \mathbb{N}$ such that for all $k\in \mathbb{N}$
\[\sup_{n\in \mathbb{N}} {\frac{\|Te_{n}\|_{k}}{\|e_{n}\|_{m}}<\infty }.\]
\end{itemize}
\end{lem}
A Fr\'echet space $E$ is Montel if each bounded set in $E$ is relatively compact. Every power series space is Montel, see Theorem 27.9 of \cite{MV}.
~\\~\\ Proposition 2.2 in \cite{N1} establishes that the continuity condition alone suffices to ensure that a linear operator defined solely on the basis elements extends to a well-defined operator. The precise statement of Proposition 2.2 from \cite{N1} is given below.

\begin{prop}\label{P0} Let $K(a_{n,k})$, $K(b_{n,k})$ be Köthe spaces and  $(a_{n})_{n\in \mathbb{N}}\in K(b_{n,k})$  be a sequence. Let us define a linear map $T:K(a_{n,k})\to K(b_{n,k})$ such as
$$Te_{n}=a_{n} \hspace{0.5in}\text{and}\hspace{0.5in} Tx=\sum^{\infty}_{n=1}x_{n}Te_{n}$$
for every $\displaystyle x=\sum^{\infty}_{n=1} x_{n}e_{n}$ and $n\in \mathbb{N}$. If the  continuity condition 
\[\forall k\in \mathbb{N} \quad\quad \exists m\in \mathbb{N}\quad\quad\quad\quad \sup_{n\in \mathbb{N}}\frac{ \|Te_{n}\|_{k}}{\|e_{n}\|_{m}}<\infty \]
holds, then $T$ is well-defined and continuous operator.
\end{prop}

Let $T$ be a continuous linear operator on a Fr\'echet space $E$. $T^{[k]}$ denotes the k-th Ces\'aro mean given by
 $$\frac{1}{k}\sum^{k}_{m=1} T^{m}.$$
\begin{defn}\label{D1}Let $E$ be a Fr\'echet space. An operator $T\in \mathcal{L}(E)$ is called 
\begin{itemize}
\item \textbf{topologizable} if for every $p\in \mathbb{N}$ there exists $q\in \mathbb{N}$ such that for every $k\in \mathbb{N}$ there is $M_{k,p}>0$ satisfying
$$\|T^{k}(x)\|_{p}\leq M_{k,p}\|x\|_{q}$$
for every $x\in E$.
\item \textbf{m-topologizable} if for every $p\in \mathbb{N}$ there exist $q\in \mathbb{N}$ and $C_{p}\geq 1$ such that 
$$\|T^{k}(x)\|_{p}\leq C^{k}_{p}\|x\|_{q}$$
holds for every $k\in \mathbb{N}$ and $x\in E$.
\item \textbf{power bounded} if for every $p\in \mathbb{N}$ there exist $q\in \mathbb{N}$ and $C_{p}\geq 1$ such that
$$\|T^{k}(x)\|_{p}\leq C_{p}\|x\|_{q}$$
holds for every $k\in \mathbb{N}$ and $x\in E$.
\item \textbf{Ces\'aro bounded }
if for every $p\in \mathbb{N}$ there exist $q\in \mathbb{N}$ and $C_{p}\geq 1$ such that 
$$\|T^{[k]}(x)\|_{p}\leq C_{p}\|x\|_{q}$$
holds for every $k\in \mathbb{N}$ and $x\in E$.
\end{itemize}
\end{defn}

The following implications are straightforward:
$$\text{power bounded} \hspace{0.1in}\Rightarrow\hspace{0.1in}\text{m-topologizable} \hspace{0.1in}\Rightarrow\hspace{0.1in}\text{topologizable}.$$

The concepts of topologizable and m-topologizable operators were defined and studied by \.{Z}elazko in \cite{Z} to establish a suitable locally convex topology on the space of all continuous linear operators $\mathcal{L}(E)$ acting on a non-normable locally convex space $E$. In \cite{B}, Bonet provides examples of non-normable locally convex spaces such that all operators defined on them are topologizable, or even m-topologizable. Furthermore, he established that when $E$ is a Fr\'echet space without a continuous norm, there is an operator on 
E, which is not topologizable. He also proved that every bounded operator $T\in \mathcal{L}(E)$ on a locally convex space $E$ is m-topologizable \cite[Proposition 7]{B}. Therefore, every compact operator is m-topologizable. As established in Theorem 3.5 and Theorem 3.13 of \cite{N2}, the Hankel operator $H_{\theta}: \Lambda_{\infty}(\alpha) \to \Lambda_{\infty}(\alpha)$ is compact for $\theta\in \Lambda_{\infty}(\alpha)$, similarly, $H_{\theta}: \Lambda_{1}(\alpha) \to \Lambda_{1}(\alpha)$ is compact for $\theta\in (\Lambda_{1}(\alpha))^{\prime}$. Hence, these operators are m-topologizable. In contrast, the compactness of Toeplitz operators defined on power series is observed only in very limited situations. Therefore, we will disregard their compactness and instead investigate whether these operators are topologizable and power-bounded.
~\\~\\
\.{Z}elazko in \cite[Proposition 7 and Proposition 12]{Z} showed that the sum and the product of two commuting topologizable (respectively, m-topologizable) operators are topologizable (respectively, m-topologizable) too. In Example 13 of \cite{Z}, \.{Z}elazko provides a counterexample showing that the sum of two m-topologizable operators, which do not commute, need not be topologizable.

If $(E,\|\cdot\|)$ is a normed space, then it is clear that every bounded operator $T: E\to E$ is m-topologizable.  This naturally leads us to consider a stronger condition suitable for Fréchet spaces, which we introduce as the notion of \textbf{strongly tameness}. With this definition, it is immediately apparent that the class of strongly tame operators is closed under addition and product.  

\begin{defn}\label{FSC1} Let $E$ be a Fr\'echet space and $T:E\to E$ be  a linear operator. We say that $T$ is \textbf{strongly tame} if for every $p\in \mathbb{N}$ there exists a constant $C_{p}>0$ satisfying
$$\|Tx\|_{p}\leq C_{p} \|x\|_{p}$$
for every $x\in E$.
\end{defn}

Every strongly tame operator is m-topologizable. If $C_{p}\leq 1$ for every $p\in\mathbb{N}$, then the strongly tame operator is power bounded. Moreover, since strongly tameness implies m-topologizability, the sum and product of strongly tame operators are themselves m-topologizable.


In the Lemma below, we demonstrate that the notions introduced in Definition \ref{D1} and Definition \ref{FSC1} can be verified for Köthe spaces by considering only the basis of Köthe spaces.

\begin{lem}\label{L1} 
Let $K(a_{n,k})$ be a Köthe space, and let $T: K(a_{n,k}) \to K(a_{n,k})$ be a linear continuous operator.
\begin{itemize}
\item[i.] $T$ is topologizable if and only if  for every $p\in \mathbb{N}$ there exists $q\in \mathbb{N}$ such that for every $k\in \mathbb{N}$ there is $M_{k,p}>0$ satisfying 
$$\|T^{k}(e_{n})\|_{p}\leq M_{k,p}\|e_{n}\|_{q}$$ 
for every $n\in \mathbb{N}$.
\item[ii.] $T$ is m-topologizable if and only if for every $p\in \mathbb{N}$ there exist  $q\in \mathbb{N}$ and $C_{p}\geq 1$ such that 
$$\|T^{k}(e_{n})\|_{p}\leq C_{p}^{k}\|e_{n}\|_{q}$$
holds for every $k,n\in  \mathbb{N}$.
\item[iii.] $T$ is power bounded if and only if  for every $p\in \mathbb{N}$ there exist $q\in \mathbb{N}$ and $C_{p}\geq 1$ such that
$$\|T^{k}(e_{n})\|_{p}\leq C_{p} \|e_{n}\|_{q}$$
holds for every $k,n\in  \mathbb{N}$.
\item[iv.] $T$ is Ces\'aro bounded if and only if  for every $p\in \mathbb{N}$ there exist $q\in \mathbb{N}$ and $C_{p}\geq 1$ such that 
$$\|T^{[k]}(e_{n})\|_{p}\leq C_{p}\|e_{n}\|_{q}$$
holds for every $k,n\in  \mathbb{N}$. 
\item[v.]  $T$ is strongly tame if and only if for every $p\in \mathbb{N}$ there exists a constant $C_{p}>0$ satisfying
$$\hspace{1.55in}\|Te_{n}\|_{p}\leq C_{p} \|e_{n}\|_{p}\hspace{1.25in}$$
holds for every $n\in \mathbb{N}$.
\end{itemize}
\end{lem}
\begin{proof} The proof of (i) is given below. The remaining parts (ii), (iii), (iv), and (v) can be established in a similar manner. The sufficient part is clear. 
Now, assume that the following condition is satisfied:  for every $p\in \mathbb{N}$, there exists $q\in \mathbb{N}$ such that for every $k\in \mathbb{N}$, there is an $M_{k,p}$ satisfying
$$\|T^{k}(e_{n})\|_{p}\leq M_{k,p}\|e_{n}\|_{q}$$ for every $n\in  \mathbb{N}$. Let $\displaystyle x=\sum^{\infty}_{n=1}x_{n}e_n\in K(a_{n,k})$. Then, it follows
\begin{equation*}
\begin{split}
\|T^{k}(x)\|_{p}&= \bigg\| \sum^{\infty}_{n=1} x_{n} T^{k}(e_{n})\bigg\|_{p} \leq \sum^{\infty}_{n=1} |x_{n}| \|T^{k}(e_{n})\|_{p} \\ &\leq M_{k,p}\sum^{\infty}_{n=1}|x_{n}|\|e_{n}\|_{q}=M_{k,p}\|x\|_{q}
\end{split}
\end{equation*}
So $T$ is topologizable.
\end{proof} 
 \begin{defn} A continuous linear operator $T$ on a locally convex Hausdorff space $E$ is called \textbf{mean ergodic} if the limits
$$Px=\lim_{n\to \infty} \frac{1}{n}\sum^{n}_{m=1} T^{m}x, \hspace{0.15in} x\in E$$
exist in $E$. If the convergence is uniform on bounded subsets of E, then T is
called \textbf{uniformly mean ergodic}.
\end{defn}

Let $E$ be a locally convex Hausdorff space. $E$ is called (uniformly) mean ergodic if every power bounded operator on $E$ is (uniformly) mean ergodic. The following result was established by Albanese, Bonet and Ricker in \cite[Proposition 2.13]{ABR}:

\begin{thm}\label{ABR} Let $K(a_{n,k})$ be a Köthe space. Then the following assertions are equivalent.
\begin{itemize}
\item[(I)] $K(a_{n,k})$ is mean ergodic.
\item[(ii)]  $K(a_{n,k})$ is uniformly mean ergodic.
\item[(iii)]  $K(a_{n,k})$ is Montel.
\end{itemize}
\end{thm}

Given that all power series spaces are Montel, any power bounded operator defined on power series spaces is (uniformly) mean ergodic.

In Theorem 2.5 of \cite{K1}, Kalmes and Santacreu established that for operators on Montel spaces, (uniform) mean ergodicity is equivalent to being Cesàro bounded and satisfying a certain orbit growth condition.

\begin{thm}[Theorem 2.5,\cite{K1}] \label{KT} Let $E$ be a Montel space and let $T\in \mathcal{L}(E)$.
\begin{itemize}
\item[(a)] $T$ is mean ergodic if and only if $T$ is uniformly mean ergodic.
\item[(b)] The following are equivalent.
\begin{itemize}
\item[(i)] $T$ is Ces\'aro bounded and $\displaystyle \lim_{n\to \infty} \frac{T^{n}}{n}=0,$ pointwise in $E$.
\item[(ii)] $T$ is mean ergodic on $E$.
\item[(iii)] $T$ is uniformly mean ergodic on $E$.
\end{itemize}
\end{itemize}
\end{thm}
\section{Toeplitz Operators on Power Series Spaces}

In [6], Jasiczak defined a Toeplitz operator on the space of all entire functions $H(\mathbb{C})$ as a continuous linear operator whose matrix with respect to the basis $(z^{n})_{n\in \mathbb{N}}$ is Toeplitz. 
Jasiczak established that such operators admit a representation of the form $\mathcal{C}M_F$, where $M_F$ denotes the multiplication operator by the symbol $F$, and $\mathcal{C}$ is a suitably defined Cauchy transform. 
This representation can be regarded as an analog of the characterization of Toeplitz operators on the Hardy space.
In this setting, the symbol $F$ of a Toeplitz operator belongs to  $\mathfrak{S}(\mathbb{C})=H(\mathbb{C})\oplus H_{0}(\infty)$ where  $H_{0}(\infty)$ denotes the space of all germs of holomorphic functions vanishing at $\infty$.  According to the Köthe–Grothendieck–da Silva characterization, $H_{0}(\infty)$ is the strong dual of $H(\mathbb{C})$, $(H(\mathbb{C}))^{\prime}_{b}$. To proceed, we take a closer look at how the operator is defined. Let $U$ be an open simply connected neighborhood of $\infty$ in the Riemann sphere $\mathbb{C}_{\infty}$, F be a holomorphic function in the punctured neighborhood $U\setminus \{\infty\}$ of the point $\infty$, $f\in H(\mathbb{C})$ and $z\in \mathbb{C}$.  By selecting a 
 $C^{\infty}$ smooth Jordan curve  $\gamma\subset (U\setminus \{\infty\})$ so that the point $z$ is in the interior $I(\gamma)$ of the curve and the connected set $\mathbb{C}\setminus U$ is contained in $I(\gamma)$, Jasiczak defined 
\begin{equation}\label{T}(T_{F}f)(z):=\frac{1}{2\pi i} \int_{\gamma} \frac{F(\xi)\cdot f(\xi)}{\xi-z}d\xi.
\end{equation}
Since for any 
$z\in \mathbb{C}$ such a curve $\gamma$ can be chosen in this way, and by Cauchy's theorem the value of the integral does not depend on the choice of $\gamma$, $(T_{F}f)(z)$ is well-defined for every $z\in \mathbb{C}$, and defines a holomorphic function. He demonstrated that the infinite matrix consisting of the Taylor coefficients of the functions $T_{F}\xi^{n}$, in successive columns, is a Toeplitz matrix.  Furthermore, he proved in Theorem 2.2 of \cite{J1}  that if a linear operator $T:H(\mathbb{C})\to H(\mathbb{C})$ whose matrix is Toeplitz is of the form (\ref{T}). The formal statement of Theorem 2.2 from \cite{J1} is given below.


\begin{thm}\label{J1} The following conditions are equivalent:
\begin{itemize}
\item[(i)] $T:H(\mathbb{C})\to H(\mathbb{C})$ is a continuous linear operator with the matrix 
\[
\begin{pmatrix}
a_{0} &a_{-1}&a_{-2}&\cdots \\
a_{1}&a_{0}&a_{-1}&\cdots\\
a_{2}&a_{1}&a_{0}&\cdots\\
\vdots&\vdots&\vdots&\ddots
\end{pmatrix}
\]
for some complex numbers $a_{n}$, $n\in\mathbb{Z}$;
\item[ii] There exists a function $F$, which is holomorphic in a punctured neighborhood $U\setminus \{\infty\}$ of $\infty$ in $\mathbb{C}_{\infty}$ such that $T=T_{F}$. Then
$$a_{n}=\frac{1}{2\pi i} \int_{\gamma} F(z)z^{-n-1}dz, \hspace{0.1in} n\in \mathbb{Z},$$
where $\gamma$ is a $C^{\infty}$ smooth Jourdan curve in $U\setminus \{\infty\}$ (the set U may be assumed to be simply connected in $\mathbb{C}_{\infty}$) such that $\mathbb{C}\setminus U\subset I(\gamma)$ and $0\in I(\gamma)$.
\end{itemize}
\end{thm}

In \cite{J2}, Jasiczak extended his approach to Toeplitz operators acting on spaces of holomorphic functions on domains in $\mathbb{C}$, the boundary of which consists of finitely many $C^\infty$ smooth Jordan curves. For instance, he provided a characterization of Toeplitz operators on the unit disk, which is formally stated in Corollary 4.1 of \cite{J2}. In the unit disc case, the symbol space $\mathfrak{S}(\mathbb{D})$ is isomorphic to $H(\mathbb{D})\oplus (H(\mathbb{D}))^{\prime}_{b}$. The formal statement of Corollary 4.1 in \cite{J2} is given below

\begin{thm}\label{J2} Let $T:H(\mathbb{D})\to H(\mathbb{D})$ be a continuous linear operator. The following conditions are equivalent:
\begin{itemize}
\item[(i)] The matrix of the operator $T$ with respect to the Schauder basis $\{z^{n}\}_{n\in \mathbb{N}_{0}}$ is a Toeplitz matrix. That is, there exists a sequence of numbers $\{a_{n}\}_{n\in \mathbb{Z}}$ such that for $m\in \mathbb{N}_{0}$,
$$T(\xi^{m})(z)=a_{-m}+a_{-m+1}z+a_{-m+2}z^{2}+\cdots$$
\item[(ii)] There exists a function $F$ holomorphic in some annulus $\{R<|z|<1\}$, $R<1$ such that for every $f\in H(\mathbb{D})$,
\begin{equation}\label{TFF}
Tf(z)=T_{F}f(z)=\frac{1}{2\pi i} \int_{|\xi|=r} \frac{F(\xi)\cdot f(\xi)}{\xi-z}d\xi.
\end{equation}
with $r<1$ such that $r>\max\{|z|,R\}$. Then for $m\in \mathbb{Z}$,
$$a_{m}=\frac{1}{2\pi i}\int_{|\xi|=r} F(\xi)\xi^{-m-1}d\xi.$$
\end{itemize}
\end{thm}

We would like to note that, as shown in Theorems \ref{J1} and \ref{J2}, the entries of the Toeplitz matrix come from the Laurent series coefficients of the symbol 
$F$ centered at $z=0$.

The space of all entire functions, $H(\mathbb{C})$, and the space of all holomorphic functions on the unit disc, $ H(\mathbb{D})$, are isomorphic power series spaces of infinite type $\Lambda_{\infty}(n)$ and of finite type $\Lambda_{1}(n)$, respectively.
In \cite{N1}, the author extended Jasiczak’s approach \cite{J1,J2} by defining Toeplitz operators on more general power series spaces of finite or infinite type. In this context, an operator $T:\Lambda_{r}(\alpha)\to \Lambda_{r}(\alpha)$, $r\in \{1,\infty\}$ is called as a Toeplitz operator if its matrix is a Toeplitz matrix with respect to the canonical basis $\{e_{n}\}_{n\in \mathbb{N}}$. 

A Toeplitz matrix
\[
\begin{pmatrix}
a_{0} &a_{-1}&a_{-2}&\cdots \\
a_{1}&a_{0}&a_{-1}&\cdots\\ a_{2}&a_{1}&a_{0}&\cdots\\
\vdots&\vdots&\vdots&\ddots
\end{pmatrix}
\]
can be written as a sum of a lower triangular Toeplitz matrix and an upper triangular Toeplitz matrix in the following way:
\[
\begin{pmatrix}
a_{0} &a_{-1}&a_{-2}&\cdots \\
a_{1}&a_{0}&a_{-1}&\cdots\\ a_{2}&a_{1}&a_{0}&\cdots\\
\vdots&\vdots&\vdots&\ddots
\end{pmatrix} =
\begin{pmatrix}
\theta_{0} &0&0&\cdots \\
\theta_{1}&\theta_{0}&0&\cdots\\ \theta_{2}&\theta_{1}&\theta_{0} &\cdots\\
\vdots&\vdots&\vdots&\ddots
\end{pmatrix}
+
\begin{pmatrix}
 \beta_{0}&\beta_{1}&\beta_{2}&\cdots \\
0&\beta_{0}&\beta_{1}&\cdots\\ 0&0&\beta_{0}&\cdots\\
\vdots&\vdots&\vdots&\ddots
\end{pmatrix}
\]
where $a_{0}=\theta_{0}+ \beta_{0}$ and
\[\theta_{n}=a_{n}\hspace{0.2in}\text{and}\hspace{0.2in} \beta_{n}=a_{-n}.\]
for every $n\in \mathbb{N}$.
This allows us to associate any Toeplitz matrix with a pair of sequences: one corresponding to its lower triangular part and the other to its upper triangular part. 
Denoting by $T_{\theta,\beta}$ the operator associated with the Toeplitz matrix determined by the sequences $\theta$ and $\beta$, one can express this operator as the sum of two operators defined by the lower triangular and the upper triangular Toeplitz matrix, respectively. That is,
\[T_{\theta,\beta}=\widehat{T}_{\theta}+\widecheck{T}_{\beta}\]
where $\widehat{T}_{\theta}$ denotes the operator associated with the lower triangular Toeplitz matrix, and $\widecheck{T}_{\beta}$ denotes the operator associated with the upper triangular Toeplitz matrix. Consistent with Jasiczak's framework,
the symbol for a Toeplitz operator $T_{\theta,\beta}$ on power series spaces $\Lambda_{r}(\alpha)$, $r\in \{1,\infty\}$ belongs to $\Lambda_{r}(\alpha)\oplus(\Lambda_{r}(\alpha))^{\prime}$, where $\theta$ and $\beta$ belong to $\Lambda_{r}(\alpha)$ and $(\Lambda_{r}(\alpha))^{\prime}$, respectively, $r\in \{1,\infty\}$.

We now proceed to examine the definitions of these operators in detail.

\subsection{The Convolution Operators $\widehat{T}_{\theta}$ Acting on Power Series Spaces}

Let $\theta=(\theta_{n})_{n\in \mathbb{N}_{0}}$ be any sequence. The lower triangular Toeplitz matrix defined by $\theta$ is
\[
\begin{pmatrix}
\theta_{0} &0&0&0&\cdots \\
\theta_{1}&\theta_{0}&0&0&\cdots\\
\theta_{2}&\theta_{1}&\theta_{0}&0&\cdots\\
\theta_{3}&\theta_{2}&\theta_{1}&\theta_{0}&\cdots\\
\vdots&\vdots&\vdots&\vdots&\ddots
\end{pmatrix}.
\]
In \cite{N1}, we considered the operator $\widehat{T}_{\theta}: K(a_{n,k})\to K(b_{n,k})$ formally defined by
$$\widehat{T}_{\theta}e_{n}=(0,\cdots,0,\theta_{0},\theta_{1},\theta_{2},\cdots)=\sum^{\infty}_{j=n}\theta_{j-n}e_{j}$$
where the image of each basis vector is taken as the $n^{th}$ column of the lower triangular Toeplitz matrix associated with $\theta$.  Therefore, for every $\displaystyle x=\sum^{\infty}_{n=1}x_{n}e_{n}\in K(a_{n,k})$, the operator  $\widehat{T}_{\theta}$ is given by 
$$ \widehat{T}_{\theta}x=\sum^{\infty}_{n=1}x_{n} \widehat{T}_{\theta}e_{n}.$$
However, the convergence of this series in $K(b_{n,k})$ is not guaranteed in general. Therefore, we analyzed conditions ensuring that the operator is well-defined and continuous, particularly in the case where one of the Köthe spaces is a power series space in \cite{N1}.

We begin with a detailed investigation of operators acting between power series spaces, as we will study the topologizability and power boundedness of these operators. 

If $\widehat{T}_{\theta}:\Lambda_{r}(\alpha)\to \Lambda_{r}(\alpha)$, $r\in \{1,\infty\}$ is well-defined, then it follows that
$$\widehat{T}_{\theta}e_{1}=\theta\in \Lambda_{r}(\alpha).$$
Hence, when working with this operator,  the sequence $\theta$ will be taken from $\Lambda_{r}(\alpha)$, $r\in \{1,\infty\}$.

The following two propositions, 3.3 and 3.4, correspond to special cases of Theorems 3.3 and 3.6 in \cite{N1}. For brevity, we omit the proofs and direct the reader to \cite{N1} for the full computations.

\begin{prop}\label{PC1} The operator $\widehat{T}_{\theta}:\Lambda_{1}(\alpha)\to \Lambda_{1}(\alpha)$ is well-defined and continuous for every $\theta\in \Lambda_{1}(\alpha)$.
\end{prop}

\begin{prop}\label{PC2} Let $\alpha=(\alpha_{n})_{n\in \mathbb{N}}$ be a stable sequence. The operator $\widehat{T}_{\theta}:\Lambda_{\infty}(\alpha)\to \Lambda_{\infty}(\alpha)$ is well-defined and continuous for every $\theta\in \Lambda_{\infty}(\alpha)$.
\end{prop}

Given that 
\begin{equation*}
\begin{split}
\widehat{T}_{\theta}x=\sum^{\infty}_{n=1}x_{n} \widehat{T}_{\theta}e_{n}=\sum^{\infty}_{n=1}\left(\sum^{n}_{j=1} x_{j}\theta_{n-j}\right)e_{n}=\sum^{\infty}_{n=1}(\theta * x)_{n}e_{n}=\theta *x
\end{split}
\end{equation*}
holds for every $\displaystyle x=\sum^{\infty}_{n=1}x_{n}e_{n}\in \Lambda_{r}(\alpha)$, $r\in\{1,\infty\}$, it is natural to call 
$\widehat{T}_{\theta}$ a \textbf{convolution operator} with symbol $\theta$, acting on power series spaces. Proposition \ref{PC1} and Proposition \ref{PC2} confirm that the operator $\widehat{T}_{\theta}$ is well-defined in power series spaces. Hence, the following statements hold.

\begin{cor}\label{CC}
For every \(\theta, \gamma \in \Lambda_1(\alpha)\), the convolution \(\theta * \gamma\) belongs to \(\Lambda_1(\alpha)\). If \(\alpha\) is stable, then for every \(\theta, \gamma \in \Lambda_\infty(\alpha)\), the convolution \(\theta * \gamma\) belongs to \(\Lambda_\infty(\alpha)\).
\end{cor}

Following the approach used for defining an operator via a lower triangular Toeplitz matrix, a similar construction can be carried out for an operator associated with an upper triangular Toeplitz matrix.
 ~\\
\subsection{The Dual Convolution Operators $\widecheck{T}_{\beta}$ Acting on Power Series Spaces}
Let $\beta=(\beta_{n})_{n\in \mathbb{N}_{0}}$ be any sequence. The upper triangular Toeplitz matrix defined by $\beta$ is
\[
\begin{pmatrix}
\beta_{0} &\beta_{1}&\beta_{2}&\beta_{3}&\cdots \\
0&\beta_{0}&\beta_{1}&\beta_{2}&\cdots\\ 0&0&\beta_{0}&\beta_{1}&\cdots\\
0&0&0&\beta_{0}&\cdots\\
\vdots&\vdots&\vdots&\vdots&\ddots
\end{pmatrix}
\]
In \cite{N1}, we considered the operator  $\widecheck{T}_{\beta}: K(a_{n,k})\to K(b_{n,k})$ formally defined by
\[\widecheck{T}_{\beta}e_{n}=(\beta_{n-1},\beta_{n-2},\dots, \beta_{1},\beta_{0},0,0,\dots)=\sum^{n}_{j=1}\beta_{n-j}e_{j}\]
where the image of each basis vector is taken as the $n^{th}$ column of the upper triangular Toeplitz matrix associated with $\beta$.  Therefore, for every $\displaystyle x=\sum^{\infty}_{n=1}x_{n}e_{n}\in K(a_{n,k})$, the operator  $\widecheck{T}_{\beta}$ is given by 
$$ \widecheck{T}_{\beta}x=\sum^{\infty}_{n=1}x_{n} \widecheck{T}_{\beta}e_{n}.$$
Again, the convergence of this series in $K(b_{n,k})$ is not guaranteed in general. Therefore, we analyzed conditions ensuring that the operator is well-defined and continuous, in the case where one of the Köthe spaces is a power series space in \cite{N1}.

Since we aim to examine topologizability and power boundedness, we begin by analyzing operators defined on power series spaces in detail.

If  $\widehat{T}_{\beta}:\Lambda_{r}(\alpha)\to \Lambda_{r}(\alpha)$, $r\in \{1,\infty\}$ is well-defined and continuous, then Proposition 4.1 of \cite{N1} implies that $\beta$ belongs to the strong dual of $\Lambda_{r}(\alpha)$, $r\in\{1,\infty\}$ .
Therefore, when working with this operator,  the sequence $\beta$ will be taken from $(\Lambda_{r}(\alpha))^{\prime}$, $r\in \{1,\infty\}$.

The following two propositions, 3.6 and 3.7, correspond to special cases of Theorems 4.6 and 4.4 in \cite{N1}. For brevity, we omit the proofs and refer the reader to \cite{N1} for the full computations.

\begin{prop}\label{PC3} Let $\alpha=(\alpha_{n})_{n\in \mathbb{N}}$ be a stable sequence, and assume that the space $\Lambda_{1}(\alpha)$ is nuclear. Then the operator $\widecheck{T}_{\beta}:\Lambda_{1}(\alpha)\to \Lambda_{1}(\alpha)$ is well-defined and continuous for every $\beta\in (\Lambda_{1}(\alpha))^{\prime}$.
\end{prop}

\begin{prop}\label{PC4}  Suppose that the power series space of infinite type $\Lambda_{\infty}(\alpha)$ is nuclear. The operator $\widecheck{T}_{\beta}:\Lambda_{\infty}(\alpha)\to \Lambda_{\infty}(\alpha)$ is well-defined and continuous for every $\beta\in (\Lambda_{\infty}(\alpha))^{\prime}$.
\end{prop}
Given that 
\begin{equation*}
\begin{split}
\widecheck{T}_{\beta}x=\sum^{\infty}_{n=1}x_{n} \widecheck{T}_{\beta}e_{n}=\sum^{\infty}_{n=1}\left(\sum^{\infty}_{j=n} x_{j}\beta_{j-n}\right)e_{ n}
\end{split}
\end{equation*}
holds for every $\displaystyle x=\sum^{\infty}_{n=1}x_{n}e_{n}\in \Lambda_{r}(\alpha)$, $r\in\{1,\infty\}$. We want to note that the coefficients 
$$\sum^{\infty}_{j=n} x_{j}\beta_{j-n} $$
differs from those in the classical convolution form as the summation starts from $j=n$ and extends to infinity. We identify 
$\widecheck{T}_{\beta}$ as a \textbf{dual convolution operator} with symbol $\beta$, acting on power series spaces. We denote this operation  by the symbol $\beta \star x$, that is,
\[(\beta\star x)_{n}=\sum^{\infty}_{j=n} x_{j}\beta_{j-n}\]
for every $n\in \mathbb{N}$.
Proposition \ref{PC3} and Proposition \ref{PC4} confirm that the operator $\widecheck{T}_{\beta}$ is well-defined in power series spaces. Hence, the following statements hold.

\begin{cor} Let $\Lambda_{r}(\alpha)$, $r\in \{1,\infty\}$ be a nuclear power series space.
\begin{itemize}
\item[(i)] $\beta \star x$ belongs to $\Lambda_{\infty}(\alpha)$ for every $x\in \Lambda_{\infty}(\alpha)$ and $\beta\in (\Lambda_{\infty}(\alpha))^{\prime}$.
\item[(ii)] $\beta \star x$ belongs to $\Lambda_{1}(\alpha)$ for every $x\in \Lambda_{1}(\alpha)$ and $\beta\in (\Lambda_{1}(\alpha))^{\prime}$ provided that $\alpha$ is stable.
\end{itemize}
\end{cor}


\subsection{Products Within Convolution and Dual Convolution Operator Classes}
In this section, we establish that the product of convolution (respectively, dual convolution) operators is again a convolution (respectively, dual convolution) operator. This is due to the fact that the matrix representation of a convolution (respectively, dual convolution) operator is a lower (respectively, upper) triangular Toeplitz matrix. This result follows from the fundamental property that the product of two lower (respectively, upper) triangular matrices is again lower (respectively, upper) triangular, and the Toeplitz structure is preserved under such matrix multiplication.

We begin by considering the product of two convolution operators.

Let $\theta,\phi\in \Lambda_{r}(\alpha)$, $r\in \{1,\infty\}$. We note that for $r=\infty$, the sequence $\alpha$ is assumed to be stable in the context of the convolution operator $\widehat{T}_{\theta}$.
We aim to show that
\begin{equation}\label{LT1}
\widehat{T}_{\phi}\widehat{T}_{\theta}=\widehat{T}_{\phi * \theta}.
\end{equation}
where $\phi*\theta$ denotes the convolution of $\beta$ and $\theta$.  We want to note that $\phi*\theta\in \Lambda_{r}(\alpha)$, $r\in\{1,\infty\}$ by Corollary \ref{CC}. Proposition \ref{PC1} and Proposition \ref{PC2} implies that $\widehat{T}_{\phi*\theta}$ is well-defined and continuous. Since all operators $\widehat{T}_{\theta}$,$\widehat{T}_{\phi}$ and $\widehat{T}_{\phi*\theta}$ are continuous, it suffices to verify that equality (\ref{LT1}) holds for the elements of basis $(e_{n})_{n\in \mathbb{N}}$. Let $n\in \mathbb{N}$ and write $\displaystyle \widehat{T}_{\theta}(e_{n})=\sum^{\infty}_{j=1}y_{j}e_{j}$ where $y_{j}=0$ if $j<n$ and $y_{j}=\theta_{j-n}$ if $j\geq n$. Hence 
\begin{equation*}
\begin{split}
\widehat{T}_{\phi}\widehat{T}_{\theta}(e_{n})&= \sum^{\infty}_{j=n}y_{j}\widehat{T}_{\phi}(e_{j}) =y_{n}\widehat{T}_{\phi}(e_{n})+\cdots+y_{m+n}\widehat{T}_{\phi}(e_{n+m})+\cdots \\ 
&= y_{n}\sum^{\infty}_{j=n}\phi_{j-n}e_{j}  +\cdots + y_{m+n}\sum^{\infty}_{j=m+n} \phi_{j-(m+n)}e_{j} +\cdots \\
&= y_{n}\phi_{0} e_{n}+\cdots + (y_{n}\phi_{m}+y_{n+1}\phi_{m-1}+\cdots+y_{n+m}\phi_{0})e_{n+m}+\cdots\\
&=\sum^{\infty}_{j=n} \gamma_{j-n}e_{j}=\widehat{T}_{\gamma}(e_{n})
\end{split}
\end{equation*}
where $\displaystyle \gamma_{m}=\sum^{m}_{i=0}y_{n+i}\phi_{m-i}=\sum^{m}_{i=0}\theta_{i}\phi_{m-i}=(\phi*\theta)_{m}$ for all $m\in \mathbb{N}_{0}$. This in turn implies that   $\widehat{T}_{\phi}\widehat{T}_{\theta}=\widehat{T}_{\phi * \theta}.$

Since $\phi*\theta=\theta*\phi$, it follows that $\widehat{T}_{\phi}\widehat{T}_{\theta}=\widehat{T}_{\theta}\widehat{T}_{\phi}$, and consequently, $\widehat{T}_{\phi}$ and $\widehat{T}_{\theta}$ commute.
\begin{prop}\label{CF} Let $\phi,\theta\in \Lambda_{r}(\alpha)$, $r\in \{1,\infty\}$. The convolution operators $\widehat{T}_{\phi}$ and $\widehat{T}_{\theta}$ commute and
$$\widehat{T}_{\phi}\widehat{T}_{\theta}=\widehat{T}_{\theta}\widehat{T}_{\phi}=\widehat{T}_{\phi * \theta}.$$
\end{prop}

Next, we investigate the product of two dual convolution operators.
Let $\beta, \psi\in (\Lambda_{r}(\alpha))^{\prime}$, $r\in \{1,\infty\}$.  We note that power series spaces are nuclear, and when $r=1$, the sequence $\alpha$ is assumed to be stable for the operators $\widecheck{T}_{\beta}$. 
We aim to show that
$$\widecheck{T}_{\beta}\widecheck{T}_{\psi}=\widecheck{T}_{ \psi*\beta}.$$
where $\beta*\psi$ denotes the convolution of $\beta$ and $\psi$. But first, we need to show that $\beta*\psi\in (\Lambda_{r}(\alpha))^{\prime}$ for $\beta, \psi\in (\Lambda_{r}(\alpha))^{\prime}$, $r\in \{1,\infty\}$.
\begin{prop} Let $\Lambda_{\infty}(\alpha)$ be a nuclear power series space of infinite type. For all $\beta, \psi \in (\Lambda_{\infty}(\alpha))^{\prime}$, their convolution $\beta * \psi$ also belongs to $(\Lambda_{\infty}(\alpha))^{\prime}$.
\end{prop}
\begin{proof} Since $\Lambda_{\infty}(\alpha)$ is nuclear, $\displaystyle \sum^{\infty}_{n=1}e^{-m_{1}\alpha_{n}}$ is convergent for some $m_{1}\in \mathbb{N}$. It follows that there exists a $D>0$ such that
\begin{equation}\label{NI}
\sum^{n}_{j=1} e^{p\alpha_{j}}\leq D e^{(p+m_{1})\alpha_{n}}.
\end{equation}
Indeed, this is evident from the following estimates:
$$\sum^{n}_{j=1} e^{p\alpha_{j}}e^{-(p+m_{1})\alpha_{n}}\leq \sum^{n}_{j=1} e^{p\alpha_{j}}e^{-(p+m_{1})\alpha_{j}}\leq \sum^{\infty}_{j=1} e^{-m_{1}\alpha_{j}} =D<+\infty.$$
Then, we get
$$
\|\widecheck{T}_{\beta} e_{n}\|_{p} \leq C_{0}e^{m_{0}\alpha_{n}} \sum^{n}_{j=1} e^{p\alpha_{j}}\leq  C_{0} D e^{p+m_{0}+m_{1}\alpha_{n}} 
=C_{0}D\|e_{n}\|_{p+m_{0}+m_{1}}
$$
for every $n\in \mathbb{N}$.

Let 
$\beta,\psi\in (\Lambda_{\infty}(\alpha))^{\prime}$. Then there exist $m_{2}, m_{3}\in \mathbb{N}$ and $C_{1}, C_{2}>0$ such that $|\beta_{n-1}|\leq C_{1}e^{m_{2}\alpha_{n}}$ and $|\psi_{n-1}|\leq C_{2}e^{m_{3}\alpha_{n}}$ for all $n\in \mathbb{N}$.  It follows that 
\begin{equation*}
\begin{split}
(\beta*\psi)_{n}&= \sum^{n}_{j=1} |\psi_{n-j}||\beta_{j-1}| \\
& \leq C_{1}C_{2}\sum^{n}_{j=1} e^{m_{2}\alpha_{n-j}}e^{m_{3}\alpha_{j-1}} \\ & \leq C_{1}C_{2}e^{m_{2}\alpha_{n}}\sum^{n}_{j=1} e^{m_{3}\alpha_{j}}\leq C_{1}C_{2}D e^{(m_{1}+m_{2}+m_{3})\alpha_{n}}
\end{split}
\end{equation*}
for all $n\in\mathbb{N}$. 
This implies that $\beta*\theta\in (\Lambda_{\infty}(\alpha))^{\prime}$.
\end{proof}
\begin{prop}\label{NP11} Let $\alpha=(\alpha_{n})_{n\in \mathbb{N}}$ be a stable sequence, and assume that the space $\Lambda_{1}(\alpha)$ is nuclear. 
For all $\beta, \psi \in (\Lambda_{1}(\alpha))^{\prime}$, their convolution $\beta * \psi$ also belongs to $(\Lambda_{1}(\alpha))^{\prime}$.
\end{prop}
\begin{proof}  Since $\Lambda_{1}(\alpha)$ is nuclear, that is, $\displaystyle \lim_{n\to\infty} \frac{\ln n}{\alpha_{n}}=0$, we obtain
\[\displaystyle  \lim_{n\to \infty} ne^{-\frac{1}{2k}\alpha_{n}}=\lim_{n\to\infty}\displaystyle  e^{-\alpha_{n}\left(\frac{\ln n}{\alpha_{n}} + \frac{1}{2k}\right)}=0\]
for every $k\in \mathbb{N}$. It follows that for every $k\in \mathbb{N}$, there exists a $D>0$ satisfying  

\begin{equation}\label{N1}
n e^{-\frac{1}{k}\alpha_{n}}\leq De^{-\frac{1}{2k}\alpha_{n}}
\end{equation}
for all $n\in \mathbb{N}$. 

Since $\alpha$ is stable, there exists a $M\in \mathbb{N}$ such that $\alpha_{2t}\leq M\alpha_{t}$ for every $t\in \mathbb{N}$. Assume that $n=2t$ or $n=2t+1$ and $1\leq j\leq n$
\[ \alpha_{n}\leq \alpha_{2t+2}\leq M\alpha_{t+1}\leq M (\alpha_{n-j+1}+\alpha_{j})\]
in this case $t+1\leq j$ or $t+1\leq n-j+1$ and then we have $\alpha_{n-j+1}+\alpha_{j}\geq \alpha_{t+1}$. Therefore the inequality 
\begin{equation}\label{S1} \alpha_{n}\leq M(\alpha_{n-j+1}+\alpha_{j})
\end{equation}
is satisfied for all $n\in \mathbb{N}$ and $j\leq n$.

Let $\beta,\theta\in (\Lambda_{1}(\alpha))^{\prime}$. Then there exist $m_{1}, m_{2}\in \mathbb{N}$ and $C_{1}, C_{2}>0$ such that $|\beta_{n-1}|\leq C_{1}e^{-\frac{1}{m_{1}}\alpha_{n}}$ and $|\psi_{n-1}|\leq C_{2}e^{-\frac{1}{m_{2}}\alpha_{n}}$ for all $n\in \mathbb{N}$.  
It follows that 
\begin{equation*}
\begin{split}
(\beta*\psi)_{n}&= \sum^{n}_{j=1} |\psi_{n-j}||\beta_{j-1}|  \leq C_{1}C_{2}\sum^{n}_{j=1}e^{-\frac{1}{m_{1}}\alpha_{n-j+1}}e^{-\frac{1}{m_{2}}\alpha_{j}} \\ & \leq C_{1}C_{2}\sum^{n}_{j=1} e^{-\frac{1}{\max\{m_{1},m_{2}\}}(\alpha_{n-j+1}+\alpha_{j})}\leq C_{1}C_{2} \sum^{n}_{j=1} e^{-\frac{1}{M\max\{m_{1},m_{2}\}}\alpha_{n}} \\
&= C_{1}C_{2} ne^{-\frac{1}{M\max\{m_{1},m_{2}\}}\alpha_{n}}\leq C_{1}C_{2}D e^{-\frac{1}{2M\max\{m_{1},m_{2}\}}\alpha_{n}}
\end{split}
\end{equation*}
for all $n\in\mathbb{N}$ which implies that $\beta*\psi\in (\Lambda_{1}(\alpha))^{\prime}$.
\end{proof}
To show that $\widecheck{T}_{\beta}\widecheck{T}_{\psi}=\widecheck{T}_{ \psi*\beta}$, it is enough to verify that $\widecheck{T}_{\beta}\widecheck{T}_{\psi}(e_{n})=\widecheck{T}_{\psi* \beta}(e_{n})$
for every $n\in \mathbb{N}$. Let $n\in \mathbb{N}$ and write  $\displaystyle \widecheck{T}_{\psi}(e_{n})=\sum^{n}_{j=1}y_{j}e_{j}$ where $y_{j}=\psi_{n-j}$ for all $1\leq j\leq n$. Then
\begin{equation*}
\begin{split}
\widecheck{T}_{\beta}\widecheck{T}_{\psi}(e_{n})&= \sum^{n}_{j=1}y_{j}\widecheck{T}_{\beta}(e_{j}) =y_{1}\beta_{0}e_{1}+ y_{2}(\beta_{1}e_{1}+\beta_{0}e_{2}) +\cdots y_{n}(\beta_{n-1}e_{1}+\cdots+\beta_{0}e_{n}) \\
&=(y_{1}\beta_{0}+y_{2}\beta_{1}+\cdots+y_{n}\beta_{n-1})e_{1}+(y_{2}\beta_{0}+y_{3}\beta_{1}+\cdots+y_{n}\beta_{n-2})e_{2}\\
&\hspace{0.15in}+\cdots +(y_{i}\beta_{0}+y_{i+1} \beta_{1}+\cdots+y_{n}\beta_{n-i})e_{i}+\cdots+y_{n}\beta_{0}e_{n}
\\
&=\bigg(\sum^{n-1}_{j=1}y_{j}\beta_{j-1}\bigg)e_{1}+\bigg(\sum^{n-2}_{j=1}y_{j+1}\beta_{j-1}\bigg)e_{2}\\
&\hspace{0.15in}+\cdots+\bigg(\sum^{n-i}_{j=1}y_{j+i}\beta_{j-1}\bigg)e_{i}+\cdots +y_{n}\beta_{0}e_{n}
\\
&=\sum^{n}_{i=1} \gamma_{n-i}e_{i}
\end{split}
\end{equation*}
where $\displaystyle \gamma_{m}=\sum^{m}_{j=1}y_{j+n-m}\beta_{j-1}=\sum^{m}_{j=1}\psi_{m-j}\beta_{j-1}=(\psi*\beta)_{m}$  for all $m\in \mathbb{N}$. It follows that $\widecheck{T}_{\beta}\widecheck{T}_{\psi}=\widecheck{T}_{\beta * \psi}$.
~\\

Since $\psi*\beta=\beta*\psi$, it follows that $\widecheck{T}_{\beta}\widecheck{T}_{\psi}=\widecheck{T}_{\psi}\widecheck{T}_{\beta}$, and consequently,  $\widecheck{T}_{\beta}$ and $\widecheck{T}_{\psi}$ commute.
\begin{prop} Let $\beta,\psi\in (\Lambda_{r}(\alpha))^{\prime}$, $r\in \{1,\infty\}$. The operators $\widecheck{T}_{\beta}$ and $\widecheck{T}_{\psi}$ commute and
$$\widecheck{T}_{\beta}\widecheck{T}_{\psi}=\widecheck{T}_{\psi}\widecheck{T}_{\beta}=\widecheck{T}_{\psi * \beta}.$$
\end{prop}

 \section{Topologizability and 
 Power Boundedness of Convolution Operators on Power Series Spaces}
In this section, we establish necessary and sufficient conditions for the convolution operators $\widehat{T}_{\theta}$ defined on power series spaces to be topologizable, m-topologizable, and power bounded.

We first examine the convolution operators $\widehat{T}_{\theta}$ defined on a power series space of finite type $\Lambda_{1}(\alpha)$.

Let $\theta\in \Lambda_{1}(\alpha)$.
For every $p,n\in \mathbb{N}$, we have
\begin{equation*}
\begin{split}
\|\widehat{T}_{\theta}e_{n}\|_{p}&= \sum^{\infty}_{j=n} |\theta_{j-n}|e^{-\frac{2}{2p}\alpha_{j}}=\sum^{\infty}_{j=n}|\theta_{j-n}|e^{-\frac{1}{2p}\alpha_{j}} e^{-\frac{1}{2p}\alpha_{j}} \\
&\leq e^{-\frac{1}{2p}\alpha_{n}}\sum^{\infty}_{j=n}|\theta_{j-n}|e^{-\frac{1}{2p}\alpha_{j}} \leq \|e_{n}\|_{2p}\|\theta\|_{2p}
\end{split}
\end{equation*}
that is,
\begin{equation}\label{FP1}
\|\widehat{T}_{\theta}e_{n}\|_{p}\leq \|\theta\|_{2p} \|e_{n}\|_{2p}
\end{equation}
holds for every $p,n\in \mathbb{N}$.
In the following, since we will consider powers of convolution operators, we use the notation  $\theta^{*k}$ 
\[
\theta^{*k} := \theta * \theta * \cdots * \theta
\]
to denote the $k$-times convolution of $\theta$ with itself. In Proposition \ref{CF}, we showed that the product of two convolution operators $\widehat{T}_{\theta}$ and $\widehat{T}_{\phi}$, is again a convolution operator $\widehat{T}_{\theta*\phi}$. By using (\ref{FP1}) and Proposition \ref{CF}, we obtain 
\begin{equation}\label{LEPF}
\|\widehat{T}^{k}_{\theta} e_{n}\|_{p} =\|\widehat{T}_{\theta^{*k}}e_{n}\|_{p}\leq \|\theta^{*k}\|_{2p}\|e_{n}\|_{2p} .
\end{equation}
for every $p,k,n\in \mathbb{N}$ and it is easy to write that
\begin{equation*}
\|\widehat{T}^{k}_{\theta} x\|_{p} =\|\widehat{T}_{\theta^{*k}}x\|_{p}\leq \|\theta^{*k}\|_{2p}\|x\|_{2p}
\end{equation*}
for every $x\in \Lambda_{1}(\alpha)$ and $p,k\in \mathbb{N}$, indeed it follows from 
\begin{equation*}
\|\widehat{T}^{k}_{\theta} x\|_{p} =\|\widehat{T}_{\theta^{*k}}x\|_{p} =\sum^{\infty}_{n=1} |x_{n}| \|\widehat{T}_{\theta^{*k}}e_{n}\|_{p} \leq \|\theta^{*k}\|_{2p} \sum^{\infty}_{n=1}|x_{n}|\|e_{n}\|_{2p}= \|\theta^{*k}\|_{2p} \|x\|_{2p}.
\end{equation*}
Moreover, we have
\begin{equation}\label{LEPT2}
\begin{split}
\|\theta^{*k}\|_{p}&=\|\widehat{T}_{\theta^{*k}}e_{1}\|_{p}=\|\widehat{T}_{\theta}(\widehat{T}_{\theta^{*(k-1)}}e_{1})\|_{p}=\\ &\leq \|\theta\|_{2p} \|\widehat{T}_{\theta^{*(k-1)}}e_{1}\|_{2p} \leq \dots\leq \prod^{k}_{i=1}\|\theta\|_{2^{i}p}
\end{split}
\end{equation}
for every $p,k\in \mathbb{N}$. Hence, we obtain that
\begin{equation}\label{E1}
\|\widehat{T}^{k}_{\theta} e_{n}\|_{p} \leq \left(\prod^{k}_{i=1}\|\theta\|_{2^{i}p}\right)\|e_{n}\|_{2p} .
\end{equation}
for every $p,k,n\in \mathbb{N}$. This ensures that $\widehat{T}_{\theta}$ is topologizable for every $\theta\in \Lambda_{1}(\alpha)$. 
\begin{thm}
The convolution operator $\widehat{T}_{\theta}:\Lambda_{1}(\alpha)\to \Lambda_{1}(\alpha)$ is topologizable for every $\theta\in \Lambda_{1}(\alpha)$. The operator  $\widehat{T}_{\theta}:\Lambda_{1}(\alpha)\to \Lambda_{1}(\alpha)$ is m-topologizable if and only if for every $p\in \mathbb{N}$, there exists a $C_{p}>0$ satisfying $$\|\theta^{*k}\|_{p}\leq C^{k}_{p}$$ for every $k\in \mathbb{N}$. 
\end{thm}
 \begin{proof} Topologizability follows immediately from the inequality (\ref{LEPF}) and Lemma \ref{L1}. If the operator  $\widehat{T}_{\theta}:\Lambda_{1}(\alpha)\to \Lambda_{1}(\alpha)$ is m-topologizable, then for every $p\in \mathbb{N}$, there exist $q\in \mathbb{N}$ and $C_{p}>0$ such that
$$\|\widehat{T}^{k}_{\theta}x\|_{p}\leq C^{k}_{p}\|x\|_{q}$$
for every $k\in \mathbb{N}$ and $x\in \Lambda_{1}(\alpha)$. This gives us that for every $p\in \mathbb{N}$, there exists a $C_{p}>0$ such that
$$\|\theta^{*k}\|_{p}=\|\widehat{T}^{k}_{\theta}e_{1}\|_{p}\leq C^{k}_{p}\|e_{1}\|_{q}\leq C^{k}_{p}$$
for every $k\in \mathbb{N}$. The converse follows from the inequality (\ref{LEPF}).
 \end{proof} 

Given that the operators $\widehat{T}_{\theta}$ and $\widehat{T}_{\phi}$ commute and their product satisfies $\widehat{T}_{\theta}\widehat{T}_{\phi}=\widehat{T}_{\theta*\phi}$, we can conclude that their product is always topologizable for every $\theta,\phi\in \Lambda_{1}(\alpha)$.

\begin{thm}\label{FP}  The convolution operator $\widehat{T}_{\theta}:\Lambda_{1}(\alpha)\to \Lambda_{1}(\alpha)$ is power bounded if and only if $\displaystyle \sup_{k\in \mathbb{N}} \|\theta^{*k}\|_{p}<\infty$ for every $p\in \mathbb{N}$.
\end{thm}
\begin{proof} If  $\displaystyle \sup_{k\in \mathbb{N}} \|\theta^{*k}\|_{p}<\infty$ for every $p\in \mathbb{N}$, the inequality (\ref{LEPF}) and Lemma \ref{L1} implies that  $\widehat{T}_{\theta}$ is power bounded. To prove the converse, we assume that 
$\widehat{T}_{\theta}$ is power bounded. Then, for every $p\in \mathbb{N}$ there exist a $q\in \mathbb{N}$ and $C_{p}\geq 1$ such that for every $k,n\in \mathbb{N}$
$$\|\widehat{T}^{k}_{\theta}(e_{n})\|_{p}\leq C_{p}\|e_{n}\|_{q}$$
holds. In particular, taking $n=1$, for every $k\in \mathbb{N}$ we get   
$$\|\theta^{*k}\|_{p}=\|\widehat{T}^{k}_{\theta}(e_{1})\|_{p}\leq C_{p} \|e_{1}\|_{q}< C_{p}$$
which shows that $\displaystyle \sup_{k\in \mathbb{N}}\|\theta^{*k}\|_{p}<\infty$ for every $p\in \mathbb{N}$. This completes the proof.
 \end{proof}
 ~\\ We would like to note that, from inequality (\ref{LEPT2}), if $\displaystyle \sup_{p\in \mathbb{N}}\|\theta\|_{p}\leq 1$ holds, then $\displaystyle \sup_{k\in \mathbb{N}} \|\theta^{*k}\|_{p}\leq 1$ for every $p\in \mathbb{N}$ and the operator $\widehat{T}_{\theta}:\Lambda_{1}(\alpha)\to \Lambda_{1}(\alpha)$ is power bounded.
 ~\\~\\
 Based on Proposition \ref{ABR} and Theorem \ref{KT}, we are able to derive the following result:

 \begin{cor}\label{CCP} Let $\theta\in \Lambda_{1}(\alpha)$. Suppose that $\displaystyle \sup_{k\in \mathbb{N}} \|\theta^{*k}\|_{p}<\infty$ for every $p\in \mathbb{N}$.
 Then  
the convolution operator $\widehat{T}_{\theta}:\Lambda_{1}(\alpha)\to \Lambda_{1}(\alpha)$ satisfies the following properties:
\begin{itemize}
    \item[1.] $\widehat{T}_{\theta}$ is power bounded.
    \item[2.] $\widehat{T}_{\theta}$ mean ergodic.
    \item[3.]  $\widehat{T}_{\theta}$ is uniformly mean ergodic.
    \item[4.] $\widehat{T}_{\theta}$ is Ces\'aro bounded and $\displaystyle \lim_{n\to \infty} \frac{\widehat{T}_{\theta}^{\hspace{0.025in}n}}{n}=0,$ pointwise in $\Lambda_{1}(\alpha)$.
\end{itemize}
 \end{cor}
 

We next turn to convolution operators on the power series spaces of infinite type to examine their properties. The stability of the sequence $\alpha$ is required to define the convolution operator $\widehat{T}_{\theta}$ on the space $\Lambda_{\infty}(\alpha)$.

Let $\theta\in \Lambda_{\infty}(\alpha)$. Since $\alpha$ is stable, there exists a constant $M\in \mathbb{N}$ such that $\alpha_{2n}\leq M\alpha_{n}$ for every $n\in \mathbb{N}.$
Since $\alpha$ is increasing, we have the following: if $j=2t\geq 2n$, 
\[ \alpha_{j}=\alpha_{2t}\leq M\alpha_{t}\leq M\alpha_{2t-n}=M\alpha_{j-n}\leq M\alpha_{j-n+1}\]
and if $j=2t+1\geq 2n$, we have $t+1\geq n$ and 
\[\alpha_{j}=\alpha_{2t+1}\leq \alpha_{2t+2}\leq M\alpha_{t+1}\leq M\alpha_{2t+2-n}=M\alpha_{j-n+1}.\]
Therefore, we have 
\[\alpha_{j}\leq M\alpha_{j-n+1}\]
for all $j\geq 2n$.
For every $p,n\in \mathbb{N}$, we have
\begin{align*}
\begin{split}
\|\widehat{T}_{\theta}e_{n}\|_{p}&=\sum^{2n-1}_{j=n} |\theta_{j-n}|e^{p\alpha_{j}}+\sum^{\infty}_{j=2n} |\theta_{j-n}|e^{p\alpha_{j}} \\ &\leq
\sum^{2n-1}_{j=n} |\theta_{j-n}|e^{p\alpha_{2n}} +\sum^{\infty}_{j=2n} |\theta_{j-n}|e^{Mp\alpha_{j-n+1}}e^{p\alpha_{j}-Mp\alpha_{j-n+1}}\\
&\leq e^{Mp\alpha_{n}}
\sum^{2n-1}_{j=n} |\theta_{j-n}|+\sum^{\infty}_{j=2n} |\theta_{j-n}|e^{Mp\alpha_{j-n+1}} \\
&\leq e^{Mp\alpha_{n}}\sum^{\infty}_{j=n} |\theta_{j-n}|e^{Mp\alpha_{j-n+1}}\\
&\leq \|\theta\|_{Mp} \|e_{n}\|_{Mp} 
\end{split}
\end{align*}
that is,
\begin{equation}\label{E22}
\|\widehat{T}_{\theta}e_{n}\|_{p}\leq \|\theta\|_{Mp} \|e_{n}\|_{Mp}
\end{equation}
holds for every $p,n\in \mathbb{N}$.
In Proposition \ref{CF}, we showed that the product of two convolution operators, say $\widehat{T}_{\theta}$ and $\widehat{T}_{\phi}$, is again a convolution operator $\widehat{T}_{\theta*\phi}$. By using (\ref{E22}), we obtain 
\begin{equation}\label{LEPI1}
\|\widehat{T}^{k}_{\theta} e_{n}\|_{p} =\|\widehat{T}_{\theta^{*k}}e_{n}\|_{p}\leq \|\theta^{*k}\|_{Mp}\|e_{n}\|_{Mp} .
\end{equation}
for every $p,k,n\in \mathbb{N}$ and it is easy to write that
\begin{equation*}
\|\widehat{T}^{k}_{\theta} x\|_{p} =\|\widehat{T}_{\theta^{*k}}x\|_{p}\leq \|\theta^{*k}\|_{Mp}\|x\|_{Mp}
\end{equation*}
for every $x\in \Lambda_{\infty}(\alpha)$ and $p,k\in \mathbb{N}$.
Moreover, we have
\begin{equation}\label{LEPI2}
\begin{split}
\|\theta^{*k}\|_{p}&=\|\widehat{T}_{\theta^{*k}}e_{1}\|_{p}=\|\widehat{T}_{\theta}(\widehat{T}_{\theta^{*(k-1)}}e_{1})\|_{p}\\ &\leq \|\theta\|_{Mp} \|\widehat{T}_{\theta^{*(k-1)}}e_{1}\|_{Mp} \leq \dots\leq  \prod^{k}_{i=1}\|\theta\|_{M^{i}p}
\end{split}
\end{equation}
for every $p,k\in \mathbb{N}$. Hence, we obtain that
\begin{equation}\label{IT1}
\|\widehat{T}^{k}_{\theta} e_{n}\|_{p} \leq \left(\prod^{k}_{i=1}\|\theta\|_{M^{i}p} \right)\|e_{n}\|_{Mp} .
\end{equation}
holds for every $p,k,n\in \mathbb{N}$. This ensures that $\widehat{T}_{\theta}$ is topologizable for every $\theta\in \Lambda_{\infty}(\alpha)$. 
\begin{thm} Let $\alpha$ be a stable sequence. 
The convolution operator $\widehat{T}_{\theta}:\Lambda_{\infty}(\alpha)\to \Lambda_{\infty}(\alpha)$ is topologizable for every $\theta\in \Lambda_{\infty}(\alpha)$.  The operator  $\widehat{T}_{\theta}:\Lambda_{\infty}(\alpha)\to \Lambda_{\infty}(\alpha)$ is m-topologizable if and only if for every $p\in \mathbb{N}$, there exists a $D_{p}>0$ satisfying $$\|\theta^{*k}\|_{p}\leq D^{k}_{p}$$ for every $k\in \mathbb{N}$. 
\end{thm}
 \begin{proof} 
 Topologizability follows immediately from the inequality (\ref{LEPI1}) and Lemma \ref{L1}. If the operator  $\widehat{T}_{\theta}:\Lambda_{\infty}(\alpha)\to \Lambda_{\infty}(\alpha)$ is m-topologizable, then for every $p\in \mathbb{N}$, there exist $q\in \mathbb{N}$ and $C_{p}>0$ such that
$$\|\widehat{T}^{k}_{\theta}x\|_{p}\leq C^{k}_{p}\|x\|_{q}$$
for every $k\in \mathbb{N}$ and $x\in \Lambda_{\infty}(\alpha)$. This gives us that for every $p\in \mathbb{N}$, we define $D_{p}:=C_{p}e^{q\alpha_{1}}$ such that
$$\|\theta^{*k}\|_{p}=\|\widehat{T}^{k}_{\theta}e_{1}\|_{p}\leq C^{k}_{p}\|e_{1}\|_{q}=C^{k}_{p}e^{q\alpha_{1}}\leq D^{k}_{p}$$
for every $k\in \mathbb{N}$. The converse follows from the inequality (\ref{LEPI1}).
 \end{proof} 
We again note that the operators $\widehat{T}_{\theta}$ and $\widehat{T}_{\phi}$ commute and their product satisfies $\widehat{T}_{\theta}\widehat{T}_{\phi}=\widehat{T}_{\theta*\phi}$, then we can conclude that their product is always topologizable for every $\theta,\phi\in \Lambda_{\infty}(\alpha)$. 
\begin{thm}\label{IP} Let $\alpha$ be a stable sequence. The convolution operator $\widehat{T}_{\theta}:\Lambda_{\infty}(\alpha)\to \Lambda_{\infty}(\alpha)$ is power bounded if and only if $\displaystyle \sup_{k\in \mathbb{N}} \|\theta^{*k}\|_{p}<\infty$ for every $p\in \mathbb{N}$.
\end{thm}
\begin{proof} If the convolution operator $\widehat{T}_{\theta}:\Lambda_{\infty}(\alpha)\to \Lambda_{\infty}(\alpha)$ is power bounded, then for every $p\in \mathbb{N}$ there exist $q\in \mathbb{N}$ and $C_{p}\geq 1$ such that
$$\|\theta^{*k}\|_{p}=\|\widehat{T}_{\theta}^{k}e_{1}\|_{p} \leq C_{p}\|e_{1}\|_{q}$$
for every $k\in \mathbb{N}$. Hence, $\displaystyle \sup_{k\in \mathbb{N}} \|\theta^{*k}\|_{p}<\infty$ for every $p\in \mathbb{N}$. For the converse, let us assume that $\displaystyle \sup_{k\in \mathbb{N}} \|\theta^{*k}\|_{p}<\infty$ for every $p\in \mathbb{N}$.
for every $k\in\mathbb{N}$. Then, the inequality (\ref{LEPI1}) and  Lemma \ref{L1} gives us that $\widehat{T}_{\theta}:\Lambda_{\infty}(\alpha)\to \Lambda_{\infty}(\alpha)$ is power bounded.
\end{proof}
 ~\\ We would like to note that, from inequality (\ref{LEPI2}), if $\displaystyle \sup_{p\in \mathbb{N}}\|\theta\|_{p}\leq 1$ holds, then $\displaystyle \sup_{k\in \mathbb{N}} \|\theta^{*k}\|_{p}\leq 1$ for every $p\in \mathbb{N}$ and the operator $\widehat{T}_{\theta}:\Lambda_{\infty}(\alpha)\to \Lambda_{\infty}(\alpha)$ is power bounded.
 ~\\~\\
We make use of the conclusions of Proposition \ref{ABR} and Theorem \ref{KT} to establish the following result.
\begin{cor}
  Let $\alpha$ be a stable sequence and $\theta\in \Lambda_{\infty}(\alpha)$. Suppose that $\displaystyle \sup_{k\in \mathbb{N}} \|\theta^{*k}\|_{p}<\infty$ for every $p\in \mathbb{N}$. Then  
the convolution operator $\widehat{T}_{\theta}:\Lambda_{\infty}(\alpha)\to \Lambda_{\infty}(\alpha)$ satisfies the following properties:
\begin{itemize}
    \item[1.] $\widehat{T}_{\theta}$ is power bounded.
    \item[2.] $\widehat{T}_{\theta}$ is  mean ergodic.
    \item[3.]  $\widehat{T}_{\theta}$ is uniformly mean ergodic.
    \item[4.] $\widehat{T}_{\theta}$ is Ces\'aro bounded and $\displaystyle \lim_{n\to \infty} \frac{\widehat{T}_{\theta}^{\hspace{0.025in}n}}{n}=0,$ pointwise in $\Lambda_{1}(\alpha)$.
\end{itemize}
 \end{cor}
\section{Topologizability and Power Boundedness  of Dual Convolution Operators on Power Series Spaces}

In this section, we establish necessary and sufficient conditions for dual convolution operators $\widecheck{T}_{\beta}$ defined on power series spaces to be topologizable, m-topologizable, and power bounded.

\begin{thm}\label{TU1} Let $\Lambda_{\infty}(\alpha)$ be a nuclear power series space and $\beta\in (\Lambda_{\infty}(\alpha))^{\prime}$. Then, the operator $\widecheck{T}_{\beta}:\Lambda_{\infty}(\alpha) \to \Lambda_{\infty}(\alpha)$ is topologizable if and only if for every $p\in \mathbb{N}$ there exists a $q\in \mathbb{N}$ so that for every $k\in \mathbb{N}$ there is a $L_{k,p}>0$ such that
$$|\beta^{*k}_{n-1}|\leq L_{k,p}e^{q\alpha_{n}}$$
for every $n\in \mathbb{N}$.
\end{thm}
\begin{proof} If the operator $\widecheck{T}_{\beta}:\Lambda_{\infty}(\alpha) \to \Lambda_{\infty}(\alpha)$ is topologizable, then for every $p\in \mathbb{N}$ there exists a 
$q\in \mathbb{N}$ so that for every $k\in \mathbb{N}$ there is a $M_{k,p}>0$ such that
$$\|\widecheck{T}^{k}_{\beta}e_{n}\|_{p}=\|\widecheck{T}_{\beta^{*k}}e_{n}\|=\sum^{n}_{j=1}|\beta^{*k}_{n-j}|e^{p\alpha_{j}}\leq M_{k,p}\|e_{n}\|_{q}$$
for every $ n\in \mathbb{N}$ and this gives us that 
$$|\beta^{*k}_{n-j}|e^{p\alpha_{j}}\leq M_{k,p}e^{q\alpha_{n}}$$
for every $1\leq j\leq n$. Taking $j=1$, we obtain
$$|\beta^{*k}_{n-1}|\leq L_{k,p} e^{q\alpha_{n}}$$
where $L_{k,p}=e^{-p\alpha_{1}}M_{k,p}$. For the converse, assume that for every $p\in \mathbb{N}$ there exists a $q\in \mathbb{N}$ so that for every $k\in \mathbb{N}$ there is a $L_{k,p}>0$ such that
$$|\beta^{*k}_{n-1}|\leq L_{k,p}e^{q\alpha_{n}}$$
holds for every $n\in \mathbb{N}$. Then we can write
\begin{equation*}
 \begin{split}
\|\widecheck{T}^{k}_{\beta}e_{n}\|_{p}&=\sum^{n}_{j=1}|\beta^{*k}_{n-j}|e^{p\alpha_{j}}=\sum^{n}_{j=1}|\beta^{*k}_{j-1}|e^{p\alpha_{n-j+1}}\\ &\leq L_{k,p}\sum^{n}_{j=1}e^{q\alpha_{j}} e^{p\alpha_{n-j+1}} \leq L_{k,p}e^{p\alpha_{n}}\sum^{n}_{j=1}e^{q\alpha_{j}} 
 \end{split}
\end{equation*}
for every $n\in \mathbb{N}$. Since $\Lambda_{\infty}(\alpha)$ is nuclear, $\displaystyle \sum^{\infty}_{n=1}e^{-m_{1}\alpha_{n}}$ is convergent for some $m_{1}\in \mathbb{N}$. It follows that there exists a $D>0$ such that
\begin{equation}\label{NI}
\sum^{n}_{j=1} e^{p\alpha_{j}}\leq D e^{(p+m_{1})\alpha_{n}}.
\end{equation}
Indeed, this is evident from the following estimates:
$$\sum^{n}_{j=1} e^{p\alpha_{j}}e^{-(p+m_{1})\alpha_{n}}\leq \sum^{n}_{j=1} e^{p\alpha_{j}}e^{-(p+m_{1})\alpha_{j}}\leq \sum^{\infty}_{j=1} e^{-m_{1}\alpha_{j}} =D<+\infty.$$
Hence for every $p\in \mathbb{N}$ there exists a $q^{\prime}=p+q+m_{1}$ such that for every $k\in \mathbb{N}$ there is $M_{k,p}=DL_{k,p}$ such that
\begin{equation*}
\begin{split}
\|\widecheck{T}^{k}_{\beta}e_{n}\|_{p} & \leq L_{k,p}e^{p\alpha_{n}}\sum^{n}_{j=1}e^{q\alpha_{j}}\leq D L_{k,p} e^{(q+p+m_{2})\alpha_{n}}\\ &= DL_{k,p} \|e_{n}\|_{p+q+m_{1}}=M_{k,p} \|e_{n}\|_{q^{\prime}}.
\end{split}
\end{equation*}
 This says that $\widecheck{T}_{\beta}$ is topologizable.
\end{proof}
The proof of the following theorem follows by analogous arguments to those in Theorem \ref{TU1} and is therefore omitted.
\begin{thm}  Let $\Lambda_{\infty}(\alpha)$ be a nuclear power series space and $\beta\in (\Lambda_{\infty}(\alpha))^{\prime}$. Then,
\begin{itemize}
\item[(i)]
$\widecheck{T}_{\beta}:\Lambda_{\infty}(\alpha) \to \Lambda_{\infty}(\alpha)$ is m-topologizable if and only if for every $p\in \mathbb{N}$ there exist $q\in \mathbb{N}$ and $D_{p}>0$ so that for every $k\in \mathbb{N}$  such that
$$|\beta^{*k}_{n-1}|\leq D^{k}_{p}e^{q\alpha_{n}}$$
for every $n\in \mathbb{N}$.
\item[(ii)] $\widecheck{T}_{\beta}:\Lambda_{\infty}(\alpha) \to \Lambda_{\infty}(\alpha)$ is power bounded if and only if for every $p\in \mathbb{N}$ there exist $q\in \mathbb{N}$ and $C_{p}\geq 1$  so that
$$|\beta^{*k}_{n-1}|\leq C_{p}e^{q\alpha_{n}}$$
for every $k,n\in \mathbb{N}$.
\end{itemize}
\end{thm}


The result below relies again on the conclusions drawn from Proposition \ref{ABR} and Theorem \ref{KT}.

 \begin{cor}
Let $\Lambda_{\infty}(\alpha)$ be a nuclear power series space and $\beta\in (\Lambda_{\infty}(\alpha))^{\prime}$. Suppose that if for every $p\in \mathbb{N}$ there exist $q\in \mathbb{N}$ and $C_{p}\geq 1$ so that
$$|\beta^{*k}_{n-1}|\leq C_{p}e^{q\alpha_{n}}$$
for every $k,n\in \mathbb{N}$.
Then, the dual convolution operator $\widecheck{T}_{\beta}:\Lambda_{\infty}(\alpha)\to \Lambda_{\infty}(\alpha)$ satisfies the following: 
\begin{itemize}
    \item[1.] $\widecheck{T}_{\beta}$ is power bounded.
    \item[2.] $\widecheck{T}_{\beta}$ is  mean ergodic.
    \item[3.]  $\widecheck{T}_{\beta}$ is uniformly mean ergodic.
    \item[4.] $\widecheck{T}_{\beta}$ is Ces\'aro bounded and $\displaystyle \lim_{n\to \infty} \frac{\widecheck{T}^{\hspace{0.025in}n}_{\beta}}{n}=0,$ pointwise in $\Lambda_{1}(\alpha)$.
\end{itemize}
 \end{cor}

We now examine the properties of the dual convolution operator defined on power series spaces of finite type.

\begin{thm}\label{TU2}  Let $\alpha=(\alpha_{n})_{n\in \mathbb{N}}$ be a stable sequence, $\beta\in (\Lambda_{1}(\alpha))^{\prime}$, and assume that the power series space of finite type $\Lambda_{1}(\alpha)$ is nuclear.  $\widecheck{T}_{\beta}:\Lambda_{1}(\alpha) \to \Lambda_{1}(\alpha)$ is topologizable if and only if for every $p\in \mathbb{N}$ there exists a $q\in \mathbb{N}$ so that for every $k\in \mathbb{N}$ there is a $L_{k,p}>0$ such that
$$|\beta^{k}_{n-1}|\leq L_{k,p}e^{-\frac{1}{q}\alpha_{n}}$$
for every $n\in \mathbb{N}$.
\end{thm}
\begin{proof} If the operator $\widecheck{T}_{\beta}:\Lambda_{1}(\alpha) \to \Lambda_{1}(\alpha)$ is topologizable, it follows that for every $p\in \mathbb{N}$ there exists a 
$q\in \mathbb{N}$ so that for every $k\in \mathbb{N}$ there is a $M_{k,p}>0$ such that
$$\|\widecheck{T}^{k}_{\beta}e_{n}\|_{p}=\sum^{n}_{j=1} |\beta^{*k}_{n-1}|e^{-\frac{1}{p}\alpha_{j}}\leq M_{k,p}\|e_{n}\|_{q}$$
for every $ n\in \mathbb{N}$, which in turn implies
$$|\beta^{*k}_{n-j}|e^{-\frac{1}{p}\alpha_{j}}\leq M_{k,p}e^{-\frac{1}{q}\alpha_{n}}$$
for all $1\leq j\leq n$ and 
$$|\beta^{*k}_{n-1}|\leq L_{k,p} e^{-\frac{1}{q}\alpha_{n}}$$
where $L_{k,p}=e^{\frac{1}{p}\alpha_{1}}M_{k,p}$. To prove the converse, suppose that for every $p\in \mathbb{N}$ there exists $q\in \mathbb{N}$ such that for every $k\in \mathbb{N}$, one can find a constant $L_{k,p}$ satisfying
$$|\beta^{*k}_{n-1}|\leq L_{k,p}e^{-\frac{1}{q}\alpha_{n}}$$
for every $n\in \mathbb{N}$. Hence, for every $n\in \mathbb{N}$, we obtain
\begin{equation*}
 \begin{split}
\|\widecheck{T}^{k}_{\beta}e_{n}\|_{p}&=\sum^{n}_{j=1}|\beta^{*k}_{n-j}|e^{-\frac{1}{p}\alpha_{j}}=\sum^{n}_{j=1}|\beta^{*k}_{j-1}|e^{-\frac{1}{p}\alpha_{n-j+1}}\\ 
 &\leq L_{k,p}\sum^{n}_{j=1}e^{-\frac{1}{q}\alpha_{j}} e^{-\frac{1}{p}\alpha_{n-j+1}} \leq  L_{k,p}\sum^{n}_{j=1}e^{-\frac{1}{\max\{q,p\}}(\alpha_{j}+\alpha_{n-j+1})}
 \end{split}
\end{equation*}
The stability of $\alpha$ ensures the existence of $M\in \mathbb{N}$ satisfying $\alpha_{n}\leq M(\alpha_{n-j+1}+\alpha_{j})$ for all $1\leq j\leq n$, $n\in \mathbb{N}$, see (\ref{S1}) in Proposition \ref{NP11} for details. Due to the nuclearity of $\Lambda_{1}(\alpha)$, for every $k\in \mathbb{N}$ there exists a $D_{k}>0$ satisfying
$$ne^{-\frac{1}{k}\alpha_{n}}\leq D_{k}e^{-\frac{1}{2k}\alpha_{n}}$$
for every $n\in \mathbb{N}$, for details, refer to (\ref{N1}) in Proposition \ref{NP11}.  Employing these properties, we now establish the following
\begin{equation*}
 \begin{split}
\|\widecheck{T}^{k}_{\beta}e_{n}\|_{p}& \leq  L_{k,p}\sum^{n}_{j=1}e^{-\frac{1}{\max\{q,p\}}(\alpha_{j}+\alpha_{n-j+1})}
\leq L_{k,p} \sum^{n}_{j=1} e^{-\frac{1}{M\max\{q,p\}}\alpha_{n}}\\ &=   L_{k,p} n e^{-\frac{1}{M\max\{q,p\}}\alpha_{n}} \leq D_{p} L_{k,p} e^{-\frac{1}{2M\max\{q,p\}}\alpha_{n}}= D_{p}L_{k,p} \|e_{n}\|_{2M\max\{q,p\}}
 \end{split}
\end{equation*}
This says that $\widecheck{T}_{\beta}$ is topologizable.
\end{proof}

The proof of the following theorem follows by analogous arguments to those in Theorem \ref{TU2} and is therefore omitted.
\begin{thm}  Let $\alpha=(\alpha_{n})_{n\in \mathbb{N}}$ be a stable sequence, $\beta\in (\Lambda_{1}(\alpha))^{\prime}$, and assume that the power series space of finite type $\Lambda_{1}(\alpha)$ is nuclear.  Then, 
\begin{itemize}
\item[(i)]
$\widecheck{T}_{\beta}:\Lambda_{1}(\alpha) \to \Lambda_{1}(\alpha)$ is m-topologizable if and only if for every $p\in \mathbb{N}$ there exist  $q\in \mathbb{N}$ and $D_{p}>0$ so that for every $k\in \mathbb{N}$ there is $L_{k,p}$ such that
$$|\beta^{*k}_{n-1}|\leq D^{k}_{p}e^{-\frac{1}{q}\alpha_{n}}$$
for every $n\in \mathbb{N}$.
\item[(ii)] $\widecheck{T}_{\beta}:\Lambda_{1}(\alpha) \to \Lambda_{1}(\alpha)$ is power bounded if and only if for every $p\in \mathbb{N}$ there exist $q\in \mathbb{N}$ and $C_{p}\geq 1$ so that 
$$|\beta^{*k}_{n-1}|\leq C_{p} e^{-\frac{1}{q}\alpha_{n}}$$
for every $k,n\in \mathbb{N}$.
\end{itemize}
\end{thm}

From Proposition \ref{ABR} and Theorem \ref{KT}, one can deduce the following:

\begin{cor}
Let $\alpha=(\alpha_{n})_{n\in \mathbb{N}}$ be a stable sequence, $\beta\in (\Lambda_{1}(\alpha))^{\prime}$, and assume that the space $\Lambda_{1}(\alpha)$ is nuclear. Suppose that
if for every $p\in \mathbb{N}$ there exist $q\in \mathbb{N}$ and $C_{p}\geq 1$ so that 
$$|\beta^{*k}_{n-1}|\leq C_{p} e^{-\frac{1}{q}\alpha_{n}}$$
for every $k,n\in \mathbb{N}$. Then the dual  convolution operator $\widecheck{T}_{\beta}:\Lambda_{1}(\alpha)\to \Lambda_{1}(\alpha)$  satisfies the following properties:
\begin{itemize}
    \item[1.] $\widecheck{T}_{\beta}$ is power bounded.
    \item[2.] $\widecheck{T}_{\beta}$ is  mean ergodic.
    \item[3.]  $\widecheck{T}_{\beta}$ is uniformly mean ergodic.
    \item[4.] $\widecheck{T}_{\beta}$ is Ces\'aro bounded and $\displaystyle \lim_{n\to \infty} \frac{\widecheck{T}^{\hspace{0.025in}n}_{\beta}}{n}=0,$ pointwise in $\Lambda_{1}(\alpha)$.
\end{itemize}
 \end{cor}

\section{Topologizability and Power Boundedness  of Toeplitz Operators $\mathbf{T_{\theta,\beta}}$ on $\mathbf{\Lambda_{1}(n)}$ and $\mathbf{\Lambda_{\infty}(n)}$}

In this section, we aim to establish necessary conditions for the Toeplitz operator $T_{\theta, \beta}$, defined on the power series spaces $\Lambda_{1}(n)$ and $\Lambda_{\infty}(n)$, to be topologizable, m-topologizable, and power bounded. This, in turn, will allow us to identify necessary conditions for the operator $T_{F}$ to be topologizable, m-topologizable, and power bounded on the spaces $H(\mathbb{C})$ and $H(\mathbb{D})$.

Since the sequence $(n)_{n\in \mathbb{N}}$
is stable and satisfies
$\displaystyle \lim_{n \to \infty}\frac{\ln n}{n}=0$, the power series spaces $\Lambda_{1}(n)$ and $\Lambda_{\infty}(n)$ are nuclear, in this setting, the operator
\begin{equation*}
\begin{split}
T_{\theta,\beta}:&\;\Lambda_{r}(n)\to \Lambda_{r}(n) \\
&T_{\theta, \beta}=\widehat{T}_{\theta}+\widecheck{T}_{\beta}
\end{split}
\end{equation*}
is well-defined and continuous for all $\theta\in \Lambda_{r}(n)$ and $\beta\in (\Lambda_{r}(n))^{\prime}$, $r\in \{1,\infty\}$. 

The operators $\widehat{T}_{\theta}$ and $\widecheck{T}_{\beta}$, and are generally non-commuting, it becomes essential to examine the behavior of $\widehat{T}_{\theta}$ and $\widecheck{T}_{\beta}$ on the spaces $\Lambda_{1}(n)$ and $\Lambda_{\infty}(n)$ in greater detail.

Definition \ref{FSC1} introduced the concept of strongly tame operators, which naturally ensures m-topologizability. Furthermore, when the constants $C_{p}$ satisfy $C_{p}\leq 1$ for all $p\in \mathbb{N}$, the operator becomes power bounded. Given that this property is preserved under both operator addition and multiplication, we aim to identify conditions guaranteeing that both 
$\widehat{T}_{\theta}$ and 
$\widecheck{T}_{\beta}$ are strongly tame. This will, in turn, ensure the desired properties for the operator $T_{\theta,\beta}=\widehat{T}_{\theta} +\widecheck{T}_{\beta}$.

\begin{prop} For every $\theta\in \Lambda_{1}(n)$, the convolution operator $\widehat{T}_{\theta}:\Lambda_{1}(n)\to \Lambda_{1}(n)$ is strongly tame since  
$$\|\widehat{T}_{\theta}e_{n}\|_{p}\leq e^{\frac{1}{2p}}\|\theta\|_{2p}\|e_{n}\|_{p}$$
holds for every $p,n\in \mathbb{N}$.
\end{prop}
\begin{proof} Let $\theta\in \Lambda_{1}(n)$. For every $p,n\in \mathbb{N}$, the inequalities 
\begin{equation*}
\begin{split}\|\widehat{T}_{\theta}e_{n}\|_{p}&=\sum^{\infty}_{j=n} |\theta_{j-n}|e^{-\frac{1}{p}j}=\sum^{\infty}_{j=n} |\theta_{j-n}|e^{-\frac{1}{2p}(j-n+1)}e^{\frac{1}{2p}(j-n+1)}e^{-\frac{1}{p}j} \\
&= \sum^{\infty}_{j=n} |\theta_{j-n}|e^{-\frac{1}{2p}(j-n+1)}e^{-\frac{1}{2p}(j+n)}e^{\frac{1}{2p}} \\
&\leq e^{\frac{1}{2p}} e^{-\frac{1}{p}n}\sum^{\infty}_{j=n}|\theta_{j-n}|e^{-\frac{1}{2p}(j-n+1)}=e^{\frac{1}{2p}}\|\theta\|_{2p}\|e_{n}\|_{p}
\end{split}
\end{equation*}
holds. This implies that $\widehat{T}_{\theta}:\Lambda_{1}(n)\to \Lambda_{1}(n)$ is strongly tame by Lemma \ref{L1}.
\end{proof}

\begin{prop} For every $\theta\in \Lambda_{\infty}(n)$, the convolution operator $\widehat{T}_{\theta}:\Lambda_{\infty}(n)\to \Lambda_{\infty}(n)$ is strongly tame since  
$$\|\widehat{T}_{\theta}e_{n}\|_{p}\leq \|\theta\|_{p}\|e_{n}\|_{p}$$
holds for every $p,n\in \mathbb{N}$.
\end{prop}
\begin{proof} Let $\theta\in \Lambda_{\infty}(n)$.  For every $p,n\in \mathbb{N}$, the inequalities
\begin{equation*}
\begin{split}\|\widehat{T}_{\theta}e_{n}\|_{p}&=\sum^{\infty}_{j=n} |\theta_{j-n}|e^{pj}=e^{pn} \sum^{\infty}_{j=n} |\theta_{j-n}|e^{p(j-n)}\\ &\leq e^{pn} \sum^{\infty}_{j=n} |\theta_{j-n}|e^{p(j-n+1)} =\|\theta\|_{p}\|e_{n}\|_{p}
\end{split}
\end{equation*}
holds. This implies that $\widehat{T}_{\theta}:\Lambda_{\infty}(n)\to \Lambda_{\infty}(n)$ is strongly tame, by Lemma \ref{L1}.
\end{proof}

\begin{prop} Let $\Lambda_{\infty}(\alpha)$ be a nuclear space and $\beta\in (\Lambda_{\infty}(\alpha_n))^{\prime}$. If $\displaystyle A=\sum^{\infty}_{n=1} |\beta_{n-1}|<\infty$, then it follows that
$$\|\widecheck{T}_{\beta}e_{n}\|_{p}\leq A\|e_{n}\|_{p}$$
for every $n,p\in \mathbb{N}$ and this implies that the dual convolution operator $\widecheck{T}_{\beta}:\Lambda_{\infty}(\alpha)\to \Lambda_{\infty}(\alpha)$ is strongly tame.
\end{prop}
\begin{proof} For every $\beta\in (\Lambda_{\infty}(\alpha))^{\prime}$ satisfying $\displaystyle A=\sum^{\infty}_{n=1} |\beta_{n-1}|<\infty$, it follows that
\begin{equation*}
\begin{split}
\|\widecheck{T}_{\beta}e_{n}\|&=\sum^{n}_{j=1}|\beta_{n-j}|\|e_{j}\|_{p}=\sum^{n}_{j=1}|\beta_{n-j}|e^{p\alpha_{j}} \leq e^{p\alpha_{n}} \sum^{n}_{j=1}|\beta_{n-j}| \\
&=e^{p\alpha_{n}} \sum^{n}_{j=1} |\beta_{j-1}|\leq A e^{p\alpha_{n}}= A \|e_{n}\|_{p}  
\end{split}
\end{equation*}
for every $n,p\in \mathbb{N}$. Hence the dual convolution operator $\widecheck{T}_{\beta}:\Lambda_{\infty}(\alpha)\to \Lambda_{\infty}(\alpha)$ is strongly tame by Lemma \ref{L1}.
\end{proof}

 
\begin{prop} Let  $\beta\in (\Lambda_{1}(n))^{\prime}$. If $\displaystyle B=\sum^{\infty}_{n=1} |\beta_{n-1}|e^{n}<\infty$, then it follows that
$$\|\widecheck{T}_{\beta}e_{n}\|_{p}\leq B\|e_{n}\|_{p}$$
for every $n,p\in \mathbb{N}$ and this implies that the dual convolution operator $\widecheck{T}_{\beta}:\Lambda_{1}(n)\to \Lambda_{1}(n)$ is strongly tame.
\end{prop}
\begin{proof}  For every $\theta\in (\Lambda_{\infty}(\alpha))^{\prime}$ satisfying $\displaystyle B=\sum_{n\in \mathbb{N}} |\beta_{n-1}|e^{n}<\infty$, it follows that
\begin{equation*}
\begin{split}
\|\widecheck{T}_{\beta}e_{n}\|_{p}&=\sum^{n}_{j=1}|\beta_{n-j}|e^{-\frac{1}{p}j}  =e^{-\frac{1}{p}n}\sum^{n}_{j=1}|\beta_{n-j}|e^{\frac{1}{p}(n-j)}
\\ &=e^{-\frac{1}{p}n} \sum^{n-1}_{j=0}|\beta_{j}|e^{\frac{1}{p}j} \leq Be^{-\frac{1}{p}n} =B\|e_{n}\|_{p}  
\end{split}
\end{equation*}
for every $n,p\in \mathbb{N}$. Hence the dual convolution operator $\widecheck{T}_{\beta}:\Lambda_{1}(n)\to \Lambda_{1}(n)$ is strongly tame by Lemma \ref{L1}.
\end{proof}

\begin{prop}\label{ST1}Let $\theta\in \Lambda_{\infty}(n)$ and $\beta\in (\Lambda_{\infty}(n))^{\prime}$. The Toeplitz operator $T_{\theta,\beta}:\Lambda_{\infty}(n)\to \Lambda_{\infty}(n)$ is 
\begin{itemize}
\item[(i)] strongly tame, and hence m-topologizable if $\displaystyle \sum^{\infty}_{n=1} |\beta_{n-1}|<\infty$.
\item[(ii)] power bounded if $\displaystyle \sup_{p\in \mathbb{N}} \|\theta\|_{p}+ \sum_{n\in \mathbb{N}} |\beta_{n-1}|\leq 1$.
\end{itemize}
\end{prop}

\begin{proof} Let $\theta\in \Lambda_{\infty}(n)$ and $\beta\in (\Lambda_{\infty}(n))^{\prime}$. The convolution operator $\widehat{T}_{\theta}$ is strongly tame and
$$\|\widehat{T}_{\theta}e_{n}\|_{p}\leq \|\theta\|_{p}\|e_{n}\|_{p}$$
holds for every $p,n\in \mathbb{N}$. If $\displaystyle A=\sum_{n\in \mathbb{N}} |\beta_{n-1}|<\infty$, the dual convolution operator $\widecheck{T}_{\beta}:\Lambda_{\infty}(n)\to \Lambda_{\infty}(n)$ is strongly tame and 
$$\|\widecheck{T}_{\beta}e_{n}\|_{p}\leq A\|e_{n}\|_{p}$$
for every $n,p\in \mathbb{N}$. Hence, we obtain
\begin{equation*}
\begin{split}
\|T_{\theta,\beta}e_{n}\|_{p}&= \|\widehat{T}_{\theta}e_{n}+\widecheck{T}_{\beta}e_{n}\|  \leq \|\widehat{T}_{\theta}e_{n}\|_{p} +\|\widecheck{T}_{\beta}e_{n}\|_{p} \\
& \leq (\|\theta\|_{p}+A)\|e_{n}\|_{p}
\end{split}
\end{equation*}
for every $n,p\in \mathbb{N}$. This implies that $T_{\theta,\beta}$ is strongly tame, and hence it is m-topologizable. Moreover, if $\displaystyle \sup_{p\in \mathbb{N}} \|\theta\|_{p}+ \sum_{n\in \mathbb{N}} |\beta_{n-1}|\leq 1$, it follows directly that $T_{\theta,\beta}$ is power bounded.
\end{proof}

An analogous result holds for the power series space of finite type $\Lambda_{1}(n)$, as stated in the following proposition.
As the proof is based once again on the property that the sum and the product of strongly tame operators remain strongly tame, it is omitted.
\begin{prop}\label{ST2}
 Let $\theta\in \Lambda_{1}(n)$ and $\beta\in (\Lambda_{1}(n))^{\prime}$. The Toeplitz operator $T_{\theta,\beta}:\Lambda_{1}(n)\to \Lambda_{1}(n)$ is 
\begin{itemize}
\item[(i)] strongly tame, and hence m-topologizable if $\displaystyle \sum_{n\in \mathbb{N}} |\beta_{n-1}|e^{n}<\infty$.
\item[(ii)] power bounded if $\displaystyle \sup_{p\in \mathbb{N}} e^{\frac{1}{2p}}\|\theta\|_{2p}+\sum_{n\in \mathbb{N}} |\beta_{n-1}|e^{n}\leq 1.$
\end{itemize}
\end{prop}

The following result follows as a direct consequence of Proposition \ref{ABR} and Theorem \ref{KT}.
\begin{cor}
Let $\theta\in \Lambda_{r}(n)$ and $\beta\in (\Lambda_{r}(n))^{\prime}$, $r\in \{1,\infty\}$. 
\begin{itemize}
\item[(a)]  The Toeplitz operator $T_{\theta,\beta}:\Lambda_{1}(n)\to \Lambda_{1}(n)$ is 
\begin{itemize}
    \item[1.] power bounded,
    \item[2.] mean ergodic,
    \item[3.] uniformly mean ergodic,
    \item[4.] Ces\'aro bounded and $\displaystyle \lim_{n\to \infty} \frac{\widecheck{T}^{\hspace{0.025in}n}_{\beta}}{n}=0,$ pointwise in $\Lambda_{1}(n)$
\end{itemize}
provided that $\displaystyle \sup_{p\in \mathbb{N}} e^{\frac{1}{2p}}\|\theta\|_{2p}+\sum_{n\in \mathbb{N}} |\beta_{n-1}|e^{n}\leq 1$.
\item[(b)] The Toeplitz operator $T_{\theta,\beta}:\Lambda_{\infty}(n)\to \Lambda_{\infty}(n)$ is 
\begin{itemize}
    \item[1.] power bounded,
    \item[2.] mean ergodic,
    \item[3.] uniformly mean ergodic,
    \item[4.] Ces\'aro bounded and $\displaystyle \lim_{n\to \infty} \frac{\widecheck{T}^{\hspace{0.025in}n}_{\beta}}{n}=0,$ pointwise in $\Lambda_{\infty}(n)$
\end{itemize}
provided that $\displaystyle \sup_{p\in \mathbb{N}} \|\theta\|_{p}+ \sum_{n\in \mathbb{N}} |\beta_{n-1}|\leq 1$.
\end{itemize}
 \end{cor}

Next, we turn our attention to Toeplitz operators defined on $H(\mathbb{C})$ and $H(\mathbb{D})$. As shown in Theorems \ref{J1} and \ref{J2}, for a symbol $F$, the entries of the Toeplitz matrix associated with the operator $T_F$, defined as in (\ref{T}) and (\ref{TFF}), arise from the Laurent expansion of $F$ at the origin.
Building upon this and Propositions \ref{ST1} and \ref{ST2}, we derive the following corollaries.

\begin{cor} Let $\displaystyle F(z)=\sum^{\infty}_{n=-\infty}a_{n}z^{n}$ be a holomorphic function in a punctured neighborhood $\infty$.  The Toeplitz operator $T_{F}:H(\mathbb{C})\to H(\mathbb{C})$ is
\begin{itemize}
\item[(i)] strongly tame, and hence m-topologizable if $\displaystyle \sum^{\infty}_{n=1} |a_{-n+1}|<\infty$.
\item[(ii)] power bounded if $\displaystyle \sup_{p\in \mathbb{N}} \left(\sum^{\infty}_{n=1}|a_{n}|e^{pn}\right)+ \sum_{n\in \mathbb{N}} |a_{-n+1}|\leq 1$.
\end{itemize}
\end{cor}

\begin{cor} Let $\displaystyle F(z)=\sum^{\infty}_{n=-\infty}a_{n}z^{n}$ be a holomorphic function in some annulus $\{R<|z|<1\}$ $R<1$.  The Toeplitz operator $T_{F}:H(\mathbb{\mathbb{D}})\to H(\mathbb{\mathbb{D}})$ is
\begin{itemize}
\item[(i)] strongly tame, and hence m-topologizable if $\displaystyle \sum^{\infty}_{n=1} |a_{-n+1}|e^{n}<\infty$.
\item[(ii)] power bounded if $\displaystyle \sup_{p\in \mathbb{N}} e^{\frac{1}{2p}}\left(\sum^{\infty}_{n=1}|a_{n}|e^{2pn}\right)+ \sum_{n\in \mathbb{N}} |a_{-n+1}|e^{n}\leq 1$.
\end{itemize}
\end{cor}

We again apply the Proposition \ref{ABR} and Theorem \ref{KT} to obtain the following.

\begin{cor}
\begin{itemize}
\item[(a)] Let $\displaystyle F(z)=\sum^{\infty}_{n=-\infty}a_{n}z^{n}$ be a holomorphic function in a punctured neighborhood $\infty$.  The Toeplitz operator $T_{F}:H(\mathbb{C})\to H(\mathbb{C})$ is
\begin{itemize}
    \item[1.] power bounded,
    \item[2.] mean ergodic,
    \item[3.] uniformly mean ergodic,
    \item[4.] Ces\'aro bounded and $\displaystyle \lim_{n\to \infty} \frac{\widecheck{T}^{\hspace{0.025in}n}_{\beta}}{n}=0,$ pointwise in $\Lambda_{1}(n)$
\end{itemize}
provided that $\displaystyle \sup_{p\in \mathbb{N}} \left(\sum^{\infty}_{n=1}|a_{n}|e^{pn}\right)+ \sum_{n\in \mathbb{N}} |a_{-n+1}|\leq 1$.
\item[(b)] Let $\displaystyle F(z)=\sum^{\infty}_{n=-\infty}a_{n}z^{n}$ be a holomorphic function in some annulus $\{R<|z|<1\}$ $R<1$.  The Toeplitz operator $T_{F}:H(\mathbb{\mathbb{D}})\to H(\mathbb{\mathbb{D}})$ is 
\begin{itemize}
    \item[1.] power bounded,
    \item[2.] mean ergodic,
    \item[3.] uniformly mean ergodic,
    \item[4.] Ces\'aro bounded and $\displaystyle \lim_{n\to \infty} \frac{\widecheck{T}^{\hspace{0.025in}n}_{\beta}}{n}=0,$ pointwise in $\Lambda_{1}(n)$
\end{itemize}
provided that $\displaystyle \sup_{p\in \mathbb{N}} e^{\frac{1}{2p}}\left(\sum^{\infty}_{n=1}|a_{n}|e^{2pn}\right)+ \sum_{n\in \mathbb{N}} |a_{-n+1}|e^{n}\leq 1$.
\end{itemize}
 \end{cor}

\section{Acknowledgements}
I would like to thank Thomas Kalmes for clarifying the confusion regarding the notion of power boundedness. I want to express my sincere gratitude to William Ross for his insightful question during my presentation at the 2024 Joint Mathematics Meeting, which brought to my attention the behavior of products of operators defined by lower and upper triangular Toeplitz matrices.

\end{document}